\documentclass[11pt]{article}

\usepackage[top=3cm,bottom=3cm,left=2.5cm,right=2.5cm]{geometry}
\usepackage{graphicx} 
\usepackage[utf8]{inputenc}
\usepackage[table,xcdraw]{xcolor} 
\usepackage{float}
\usepackage{cite}
\usepackage{graphicx}
\usepackage{textcomp}
\usepackage{xcolor}
\usepackage{tikz}
\usepackage{pgfplots}
\usepackage[export]{adjustbox}
\usetikzlibrary{spy}
\usepackage{graphicx}
\usepackage{tabularx}
\usepackage{xcolor}
\usepackage{amssymb,amsmath,amsfonts}
\usepackage{bm}
\usepackage{tikz}
\usepackage{multirow}
\usepackage{booktabs}
\usepackage{blindtext}
\usepackage{multicol}
\usepackage{tikzscale}
\usepackage{numprint}
\npdecimalsign{.}
\nprounddigits{3} 
\usepackage{caption}
\usepackage{subcaption}
\usepackage{enumitem}
\usepackage{float}
\usepackage{pict2e}
\usepackage{siunitx}
\usepackage[english]{babel}
\usepackage{authblk}
\usepackage{orcidlink}
\usepackage{lipsum}
\usepackage{wrapfig}
\usepackage{array}  
\renewcommand{\arraystretch}{2.5}
\usepackage[most]{tcolorbox}
\definecolor{darkspringgreen}{rgb}{0.09, 0.45, 0.27}

\usepackage{mathtools}
\usepackage[percent]{overpic}
\usepackage{algorithm,algcompatible,amsmath}

\DeclareMathOperator*{\argmin}{arg\,min}
\algnewcommand\INPUT{\item[\textbf{Input:}]}%
\algnewcommand\PARAMETER{\item[\textbf{Parameters:}]}%
\algnewcommand\OUTPUT{\item[\textbf{Output:}]}%

\usepackage{comment}
\usepackage{array}

\usepackage{algorithm}
\usepackage{algpseudocode}
\algnewcommand{\Inputs}[1]{%
  \State \textbf{Inputs:}
  \Statex \hspace*{\algorithmicindent}\parbox[t]{.8\linewidth}{\raggedright #1}
}
\algnewcommand{\Initialize}[1]{%
  \State \textbf{Initialize:}
  \Statex \hspace*{\algorithmicindent}\parbox[t]{.8\linewidth}{\raggedright #1}
}

\usepackage[most]{tcolorbox}
\definecolor{darkspringgreen}{rgb}{0.09, 0.45, 0.27}

\newcommand{\bu}{\bm{u}}
\newcommand{\bx}{\bm{x}}
\newcommand{\bh}{\bm{h}}
\newcommand{\bz}{\bm{z}}
\newcommand{\by}{\bm{y}}

\newcommand{\Rd}{\mathbb{R}}
\newcommand{\bA}{\bm{A}}
\newcommand{\bb}{\bm{b}}
\newcommand{\dmat}{\bm{D}}

\usepackage{comment}

\makeatletter
\renewcommand\AB@affilsepx{, \protect\Affilfont}
\makeatother

\providecommand{\keywords}[1]
{
  \small	
  \textbf{\textit{Keywords---}} #1
}

\begin{document}
\title{
Why do we regularise in every iteration \\for imaging inverse problems?
}
\author[1]{Evangelos Papoutsellis\orcidlink{0000-0002-1820-9916}\thanks{corresponding author: epapoutsellis@gmail.com}} 
\author[2]{Zeljko Kereta\orcidlink{0000-0003-2805-0037}} 
\author[3]{Kostas Papafitsoros\orcidlink{0000-0001-9691-4576}}
\affil[1]{Finden Ltd}\affil[2]{University College London}\affil[3]{Queen Mary University of London}
\date{}
\maketitle 
\begin{abstract}
Regularisation is commonly used in iterative methods for  solving imaging inverse problems.
Many algorithms involve the evaluation of the proximal operator of the regularisation term in every iteration, leading to a significant computational overhead 
since such evaluation can be costly.  
In this context, the ProxSkip algorithm, recently proposed for federated learning purposes,
emerges as an  solution. It randomly skips regularisation steps, reducing the computational time of an iterative algorithm without affecting its convergence. 
Here we explore for the first time the efficacy of ProxSkip to  a variety of imaging inverse problems and we also propose a novel PDHGSkip version. Extensive numerical results highlight the potential of these methods to accelerate computations while maintaining high-quality reconstructions.\\[2pt]
\end{abstract}

\section{Introduction} \label{sec:introduction}

Inverse problems involve the process of estimating an unknown quantity $\bu^{\dagger}\in\mathbb{X}$ from indirect and often noisy measurements $\bb\in\mathbb{Y}$ obeying  $\bb = \bA \bu^{\dagger} + \bm{\eta}$. Here $\mathbb{X}, \mathbb{Y}$ denote finite dimensional spaces,
$\bA:\mathbb{X}\rightarrow\mathbb{Y}$ is a linear forward operator, $\bu^{\dagger}$ is the ground truth and $\bm{\eta}$ is a random noise component. Given $\bb$ and $\bA$, the goal is to compute an approximation $\bu$ of $\bu^{\dagger}$. Since inverse problems are typically ill-posed, prior information about $\bu$  has to be incorporated in the form of regularisation. The solution to the inverse problem is then acquired by solving 
\begin{equation}
\argmin_{\bu\in\mathbb{X}} \mathcal{D}(\bA\bu, \bb) + \alpha\mathcal{R}(\bu).
\label{eq:general_optimisation_problem}
\end{equation}
Here $\mathcal{D}$ denotes the fidelity term, 
measuring the distance between $\bb$ and the solution $\bx$ under the operator $\bA$. Regularisation term $\mathcal{R}$ promotes properties such as smoothness, sparsity, edge preservation, and low-rankness of the solution, and is weighted by a  parameter $\alpha>0$.
Classical examples for $\mathcal{R}$ in imaging include the well-known Total Variation (TV), high order extensions \cite{Bredies2020}, namely the Total Generalized Variation (TGV) \cite{Bredies2010}, Total Nuclear Variation \cite{Holt2014} and more general tensor based structure regularisation, \cite{Lefkimmiatis2013}.

In order to obtain a solution for \eqref{eq:general_optimisation_problem}, one employs iterative  algorithms such as Gradient Descent (GD) for smooth objectives or Forward-Backward Splitting (FBS) \cite{Combettes2005} for non-smooth ones. Moreover, under the general framework
\begin{align}
&\min_{\bx\in\mathbb{X}} f(\bx) + g(\bx),
\label{eq:fbs_objective}
\end{align}
the Proximal Gradient Descent (PGD) algorithm, also known as Iterative Shrinkage Thresholding Algorithm (ISTA) and its accelerated version FISTA \cite{BeckTeboulle} are commonly used when $f$ is a convex, L-smooth and $g$ proper convex. Saddle-point methods such as the Primal Dual Hybrid Gradient (PDHG) \cite{ChambollePock2011} are commonly used for non-smooth $f$.

A common property of most of these methods, see Algorithms \ref{alg:GD}--\ref{alg:FISTA} below, is the evaluation of proximal operators  related to the regulariser in every iteration, which for $\tau>0$ is defined as 
\begin{equation}
\mathrm{prox}_{\tau\mathcal{R}}(\bx):=\argmin_{\bz\in\mathbb{X}}\bigg\{\frac{1}{2}\|\bz-\bx\|_{2}^{2} + \tau\mathcal{R}(\bz)\bigg\}. 
\label{eq:proximal_definition}
\end{equation}
This proximal operator can have either a closed form solution, e.g., when $\mathcal{R}(\cdot)=\|\cdot\|_{1}$ or requires an inner iterative solver  e.g., when $\mathcal{R}(\bu)=\mathrm{TV}(\bu)$.

\begin{figure}[t!]
\noindent
\begin{minipage}[t]{0.49\textwidth}
\vspace{-0.9cm}
\begin{algorithm}[H]
\begin{algorithmic}[1]
\State {\bf Parameters:} $\gamma>0$
\State {\bf Initialize:} $\bx_0 \in \mathbb{X}$
\For{$k=0,\ldots, K-1$}
\State $\bx_{k+1} = \bx_{k} - \gamma \nabla f(\bx_{k})$ $\phantom{\mathrm{prox}_{\gamma g}}$
\EndFor  
\end{algorithmic}
\caption{GD \phantom{PGD/ISTA/FBS}}
\label{alg:GD}
\end{algorithm}
\end{minipage}\hfill
\begin{minipage}[t]{0.49\textwidth}
\centering
\vspace{-0.9cm}
\begin{algorithm}[H]
\begin{algorithmic}[1]
\State {\bf Parameters:} $\gamma >0$
\State {\bf Initialize:} $\bx_0 \in \mathbb{X}$
\For{$k=0,\ldots, K-1$}
\State $\bx_{k+1} = \mathrm{prox}_{\gamma g}(\bx_{k} - \gamma \nabla f(\bx_{k}))$
\EndFor  
\end{algorithmic}
\caption{PGD/ISTA/FBS}
\label{alg:PGD}
\end{algorithm}
\end{minipage}

\noindent
\begin{minipage}[t]{0.49\textwidth}
\centering
\begin{algorithm}[H]
\begin{algorithmic}[1]
\State {\bf Parameters:} $\gamma >0$, $t_0=1$ 
\State {\bf Initialize:} $\bx_0 \in \mathbb{X}$
\For{$k=0,\ldots, K-1$}
\State $t_{k+1}=\frac{1+\sqrt{1+4t_k^2}}{2}$, $a_k=\frac{t_{k-1}}{t_{k+1}}$
\State $\Bar{\bx}_{k+1} =\bx_k+a_k\left(\bx_k-\bx_{k+1}\right)$
\State $\bx_{k+1} = \mathrm{prox}_{\gamma g}(\Bar{\bx}_{k+1} - \gamma \nabla f(\Bar{\bx}_{k+1}))$
\EndFor  
\end{algorithmic}
\caption{FISTA}
\label{alg:FISTA}
\end{algorithm}
\end{minipage}\hfill
\begin{minipage}[t]{0.49\textwidth}
\begin{algorithm}[H]
\begin{algorithmic}[1]
\State {\bf Parameters:} $\gamma > 0$, probability $p>0$
\State {\bf Initialize:} $\bx_0, \bh_{0} \in \mathbb{X}$ 
\For{$k=0,1,\dotsc,K-1$}
\State $\hat \bx_{k+1} = \bx_k - \gamma (\nabla f (\bx_k) - {\bh_k})$ \hfill 
\If{$\mathrm{Prob}(\theta_{k} = 1) = p$} 
\State  $\bx_{k+1} = \mathrm{prox}_{\frac{\gamma}{p}g}\bigg(\hat \bx_{k+1} - \frac{\gamma}{p}{\bh_k} \bigg)$ 
\Else
$\;\;\bx_{k+1} = \hat \bx_{k+1}$ \hfill 
\EndIf
\State ${ \bh_{k+1}} = {\bh_k} + \frac{p}{\gamma}(\bx_{k+1} - \hat \bx_{k+1})$ 
\EndFor  
\end{algorithmic}
\caption{ProxSkip}
\label{alg:ProxSkip}
\end{algorithm}
\end{minipage}
\end{figure}

\subsection*{The ProxSkip algorithm}
The possibility to skip the computation of the proximal operator in some iterations according to a probability $p$, accelerating the algorithms, without affecting convergence,
was  discussed in\cite{mishchenko2022proxskip}.  There, the ProxSkip algorithm was introduced to tackle  federated learning applications which also rely on computations of  expensive proximal operators.
ProxSkip introduces a control variable $\bh_{k}$, see Algorithm \ref{alg:ProxSkip}. When the proximal step is not applied, the control variable remains constant. Hence, if at iteration $k$, no proximal step has been applied previously, the accumulated error is passed to $\bx_{k+1}$ without incurring an additional computational cost. If at the iteration $k$ the proximal step is applied, the error is reduced and the control variable will be updated accordingly.

In \cite{mishchenko2022proxskip} it was shown that the ProxSkip converges  provided that $f$ in \eqref{eq:fbs_objective} is  L-smooth and $\mu$-strongly convex, and probability $p$ satisfies
\begin{equation}
    p \geq \sqrt{\mu/L}.
    \label{eq:optimal_prob}
\end{equation}
In the case of  equality in \eqref{eq:optimal_prob}, the algorithm converges (in expectation) at a linear rate with $\gamma=\frac{1}{L}$ and the iteration complexity is $\mathcal{O}(\frac{L}{\mu}\log(\frac{1}{\varepsilon}))$. In addition, the total number of proximal evaluations (in expectation) are only  $\mathcal{O}(\frac{1}{\sqrt{p}}\log(\frac{1}{\varepsilon}))$.

\subsection*{Our contribution}
Our aim is to showcase for the first time via extended numerical experiments the computational benefits of  ProxSkip for a variety of  imaging inverse problems, including challenging real-world tomographic applications. In particular, we show that ProxSkip can outperform the accelerated version of its non-skip analogue, namely FISTA. 
At the same time, we introduce a novel PDHGSkip version of the PDHG, Algorithm \ref{alg:PDHGSkip}, which we motivate via numerical experiments. We anticipate that this will spark further research around developing  skip-versions of a variety of proximal based algorithms used nowadays. 

For all our imaging experiments we consider the following optimisation problem that contains a quadratic distance term as the fidelity term,  with the (isotropic) total variation as the regulariser, i.e.,
\begin{equation*}
    \mathrm{TV}(\bu) = 
\|\dmat\bu\|_{2,1} =  \sum|(\dmat_{y}\bu, \dmat_{x}\bu)|_{2}=  \sum \sqrt{ ((\dmat_{y}\bu)^{2} + (\dmat_{x}\bu)^{2}},
\end{equation*}
\begin{equation}
    \min_{\bu\in\mathbb{X}}  \frac{1}{2}\|\bA\bu - \bb\|_{2}^{2}+\alpha \mathrm{TV}(\bu),
 \label{mainTV}   
\end{equation}
where $\dmat\bu=(\dmat_y, \dmat_x)$ is the finite difference operator under Neumann boundary conditions.

\section{ProxSkip in imaging problems with light proximals}
\subsection{Dual TV denoising}\label{sec:dualTV}

To showcase the algorithmic properties we consider a toy example, with the dual formulation of the classical TV denoising (ROF) which reads
\begin{equation}
    \min_{\|\bm{q}\|_{2,\infty}\leq\alpha} \bigg\{\mathcal{F}(\bm{q}) := \frac{1}{2}\|\bm{\mathrm{div}}\bm{q} + \bb\|_{2}^{2} + \frac{1}{2}\|\bb\|^{2}\bigg\}.\quad\quad\mbox{(Dual-ROF).}
    \label{eq:dual_TV}
\end{equation}
where $\bm{\mathrm{div}}$ is the discrete  divergence operator such that $\bm{\mathrm{div}} = - \dmat^{T}$. The solutions $\bu^{*}$ and  $\bm{q^{*}}$   of the primal and dual (ROF) problems are linked via $\bu^{*} = \bb + \bm{\mathrm{div}}\bm{q^{*}}$. 
A simple algorithm
 to solve \eqref{eq:dual_TV} was introduced in   \cite{Chambolle2004}, based on a  
  Projected Gradient Descent (ProjGD) iteration  $\bm{q_{k+1}} = \mathcal{P}_{C}(\bm{q_{k}} + \gamma\dmat(\bm{\mathrm{div}\bm{q_{k}}}+\bb))$ which is globally convergent under a fixed stepsize $\gamma\le\frac{2}{\|\dmat\|^2}$ \cite{Aujol2009}, 
  with 
\begin{equation}
\mathcal{P}_{C}(\bx) = \frac{\bx}{\max\{\alpha, \|\bx\|_{2}\}}, \quad C=\left\{\|\bm{q}\|_{2,\infty}\leq\alpha: \bm{q}\in\Rd^{2\times d}\right\},
\label{eq:projection_c}
\end{equation}
and $d$ is the image dimension. This approach became quite popular in the following years for both its simplicity and efficiency, \cite{Pock2008}, and for avoiding computing smooth approximations of TV. The projection $\mathcal{P}_{C}$ can be identified as the proximal operator of indicator function of the feasibility set $C$. Thus, ProjGD is a special case of PGD and we can apply the ProxSkip  Algorithm \ref{alg:ProxSkip}. Note that due to the divergence operator this problem  is not strongly convex. In fact, this is the case for the majority of the problems of the type \eqref{eq:general_optimisation_problem}, typically due to the non-injectivity of $\bA$. Thus, this example also shows that the strong convexity assumption could potentially be relaxed for imaging inverse problems. 

\begin{figure}[t!]
    \centering
    \begin{subfigure}[b]{0.24\textwidth}
        \centering
        \includegraphics[width=\textwidth]{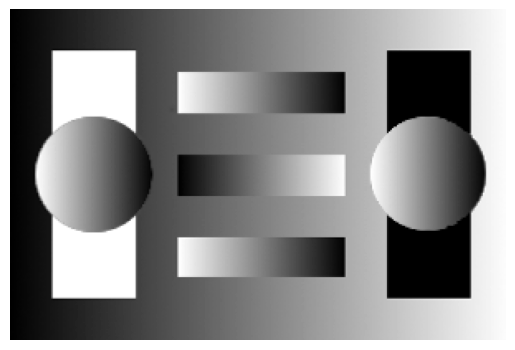}
        \label{fig_shapes_images:sub1}
    \end{subfigure}
    \begin{subfigure}[b]{0.24\textwidth}
        \centering
        \includegraphics[width=\textwidth]{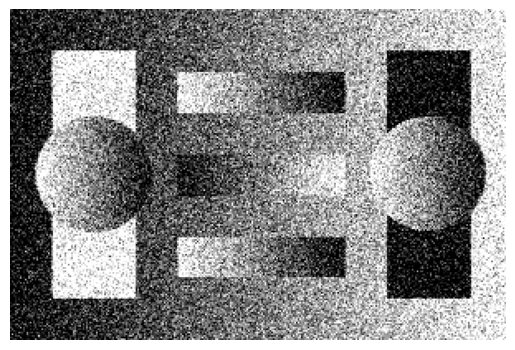}
        \label{fig_shapes_images:sub2}
    \end{subfigure}
    \begin{subfigure}[b]{0.24\textwidth}
        \centering
        \includegraphics[width=\textwidth]{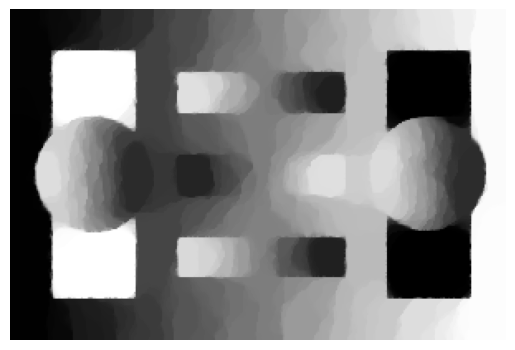}
        \label{fig_shapes_images:sub3}
    \end{subfigure}
    \begin{subfigure}[b]{0.24\textwidth}
        \centering
        \includegraphics[width=\textwidth]{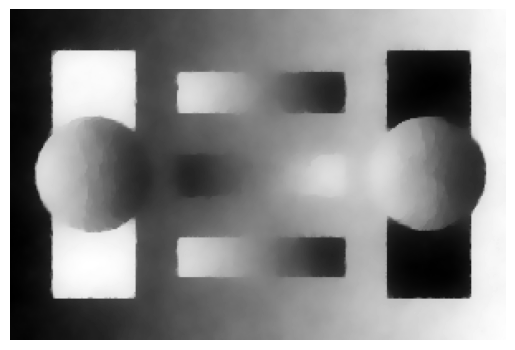}     
        \label{fig_shapes_images:sub3}
    \end{subfigure}    
    \caption{Left to right: Ground truth $\bu^{\dagger}\in\Rd^{200\times300}$.  Noisy image $\bb$, $\sigma=0.05$. Dual-ROF $u^{*}$ 
     with $\alpha=0.5$. 
    Dual-Huber-ROF $u^{*}$  (see Section \ref{sec:huber})  with $\alpha=0.55$, $\varepsilon=0.1$.
    The  parameters $\alpha$ are optimised with respect to SSIM. 
    }
    \label{fig_shapes_images:main}
\end{figure}

To ensure that any biases from algorithms under evaluation are avoided, the ``exact'' solution $\bu^{*}$ is calculated using an independent high-precision solver, in particular, the MOSEK solver from the CVXpy library \cite{mosek, diamond2016cvxpy}, see Figure \ref{fig_shapes_images:main}.
Both the ProjGD and ProxSkip algorithms use the stepsize $\gamma = \frac{1}{L} = \frac{1}{8}$, where $L$ is the Lipschitz constant of $\mathcal{F}'(\cdot)$. For every iteration, we monitor the $\ell^{2}$ error $\|\bu_k-\bu^*\|_2$ between the iterate $\bm{u}_{k} = \bb + \bm{\mathrm{div}\bm{q_{k}}}$ and estimated exact solution $\bu^{*}$. We use $p=[0.01, 0.1, 0.3, 0.5]$ and 50000 iterations as a stopping criterion.

\begin{figure}[h!]
    \centering
    \begin{subfigure}[b]{8cm}
        \centering
        \includegraphics[width=8cm]{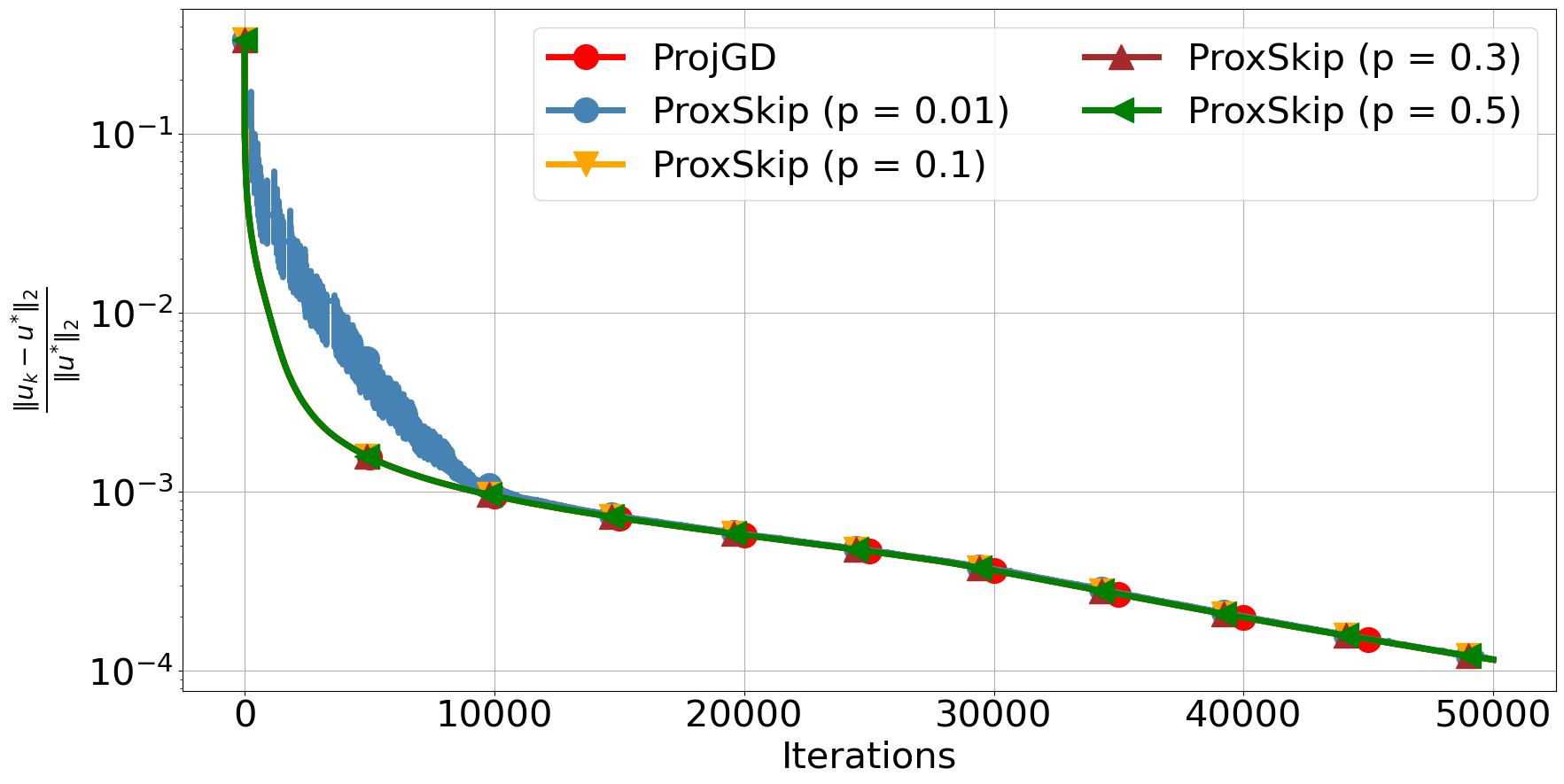}
        \label{fig_multiple_p:sub1}
    \end{subfigure}
    \begin{subfigure}[b]{8cm}
        \centering
        \includegraphics[width=8cm]{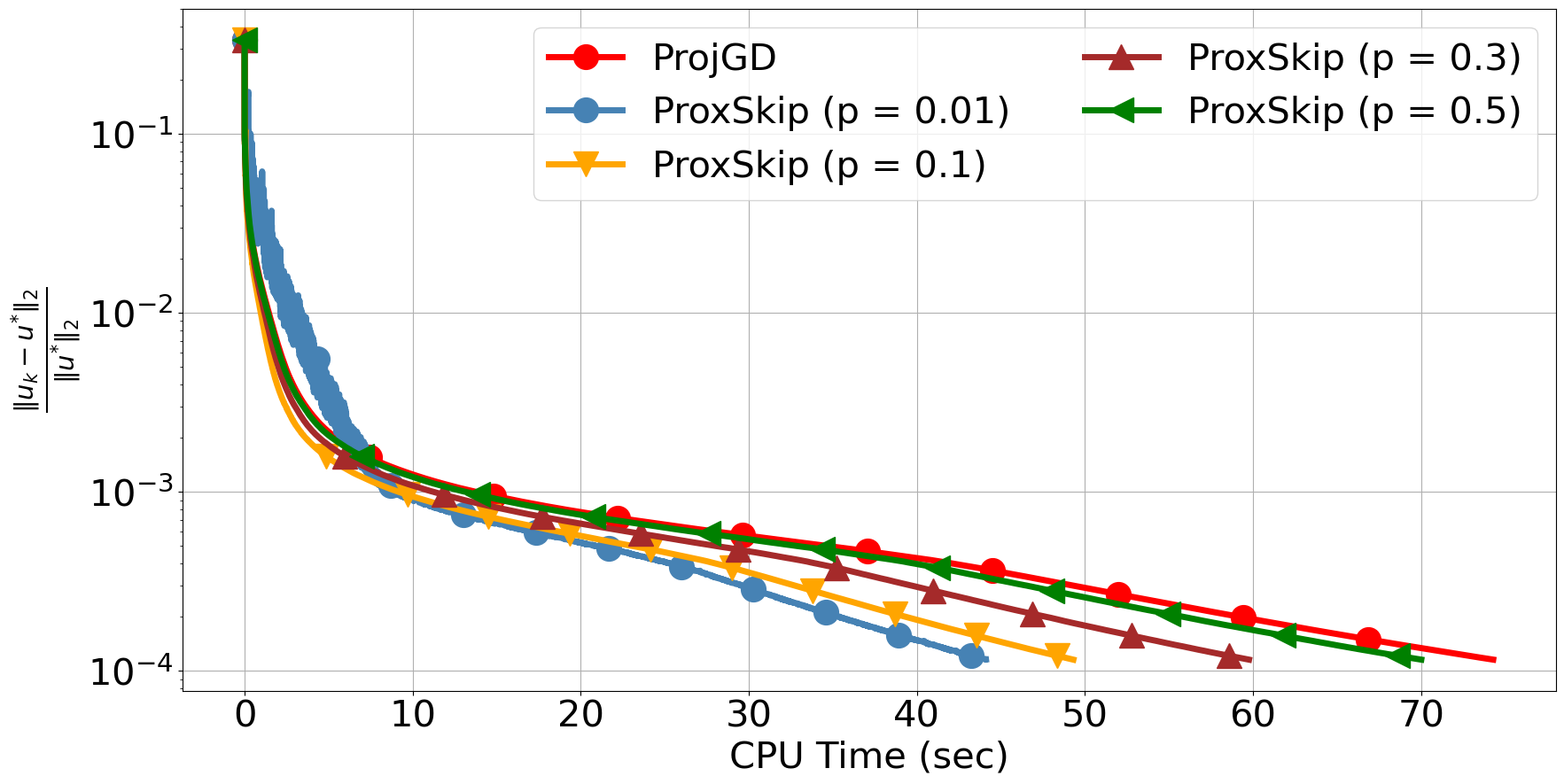}
        \label{fig_multiple_p:sub2}
    \end{subfigure}
    
    \begin{subfigure}[b]{8cm}
        \centering
        \includegraphics[width=8cm]{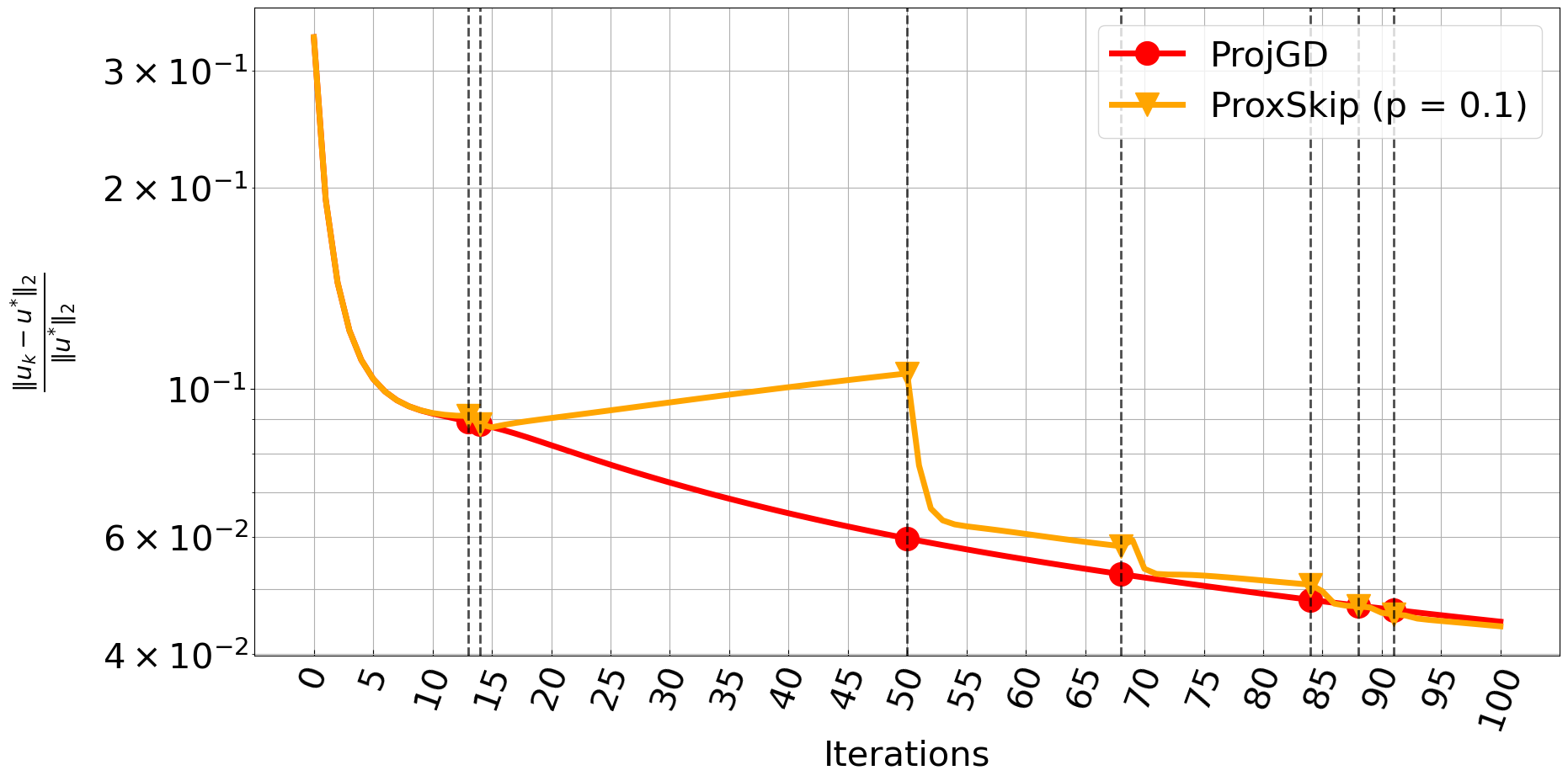}
        \label{fig_plots_iterates_objectives_denoising_zoom:sub1}
    \end{subfigure}     
    \begin{subfigure}[b]{8cm}
        \centering
        \includegraphics[width=8cm]{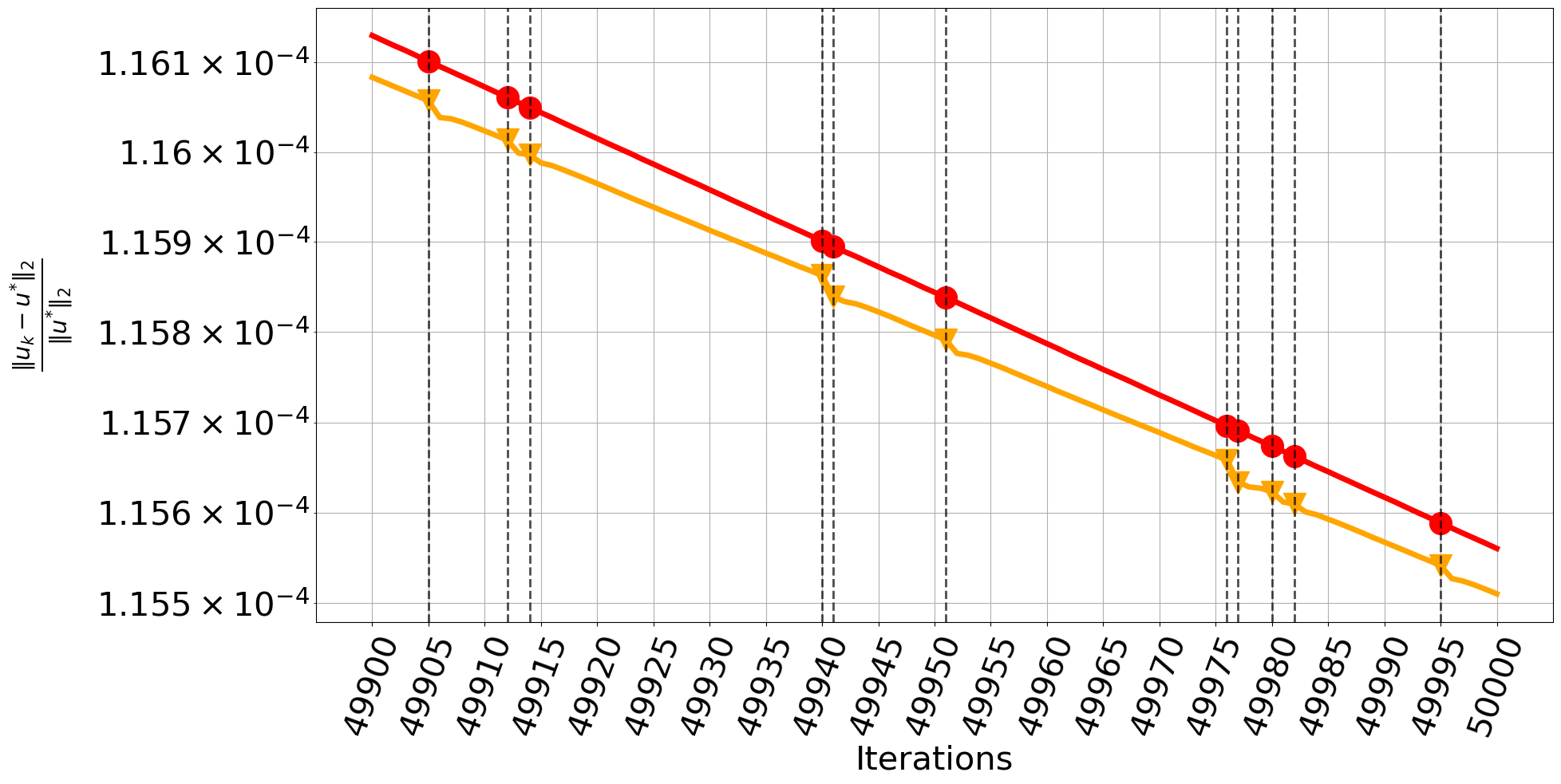}
        \label{fig_plots_iterates_objectives_denoising_zoom:sub2}
    \end{subfigure}   
    \caption{Top: Comparison of ProjGD and ProxSkip for multiple values of $p$ for \eqref{eq:dual_TV} with respect to iterations (left) and CPU time (right). Bottom: Detailed versions  for the first 100 (left) and the last 100 iterations (right) when $p=0.1$. The vertical dotted lines indicate the iterations where $\mathcal{P}_{C}(\cdot)$ is applied.}
    \label{fig_subopts_multiple_p:main}
\end{figure}

In Figure \ref{fig_subopts_multiple_p:main} (top-right), it can be observed that these two algorithms are almost identical in terms of the $\ell_2$ error. Note that ProxSkip and ProjGD coincide only when $p=1$. Indeed, one can detect some discrepancies during first 100 iterations, which quickly dissipate with only a few applications of the projection $\mathcal{P}_{C}$,
see bottom row of Figure \ref{fig_subopts_multiple_p:main}. In Figure \ref{fig_subopts_multiple_p:main} (top-left), we plot the $\ell_2$ error with respect to CPU time. The shown CPU time is the average over 30 independent runs of each algorithm. We observe a clearly superior performance of ProxSkip, for all values of $p$.
%
This serves as a first demonstration of  
the advantage of ProxSkip in terms of computational time without affecting the quality of the image. In Section \ref{sec:heavy}, we present a more emphatic computational impact using heavier proximal steps.

\subsection{Dual TV denoising with strong convexity}\label{sec:huber}

In order to be consistent with the convergence theory of ProxSkip where strong convexity is a requirement, one can add a small quadratic term to the objective function. This is a commonly used in imaging applications and allows the use of accelerated versions of first-order methods. For the \eqref{eq:dual_TV} problem this results in
\begin{equation}    \min_{\|\bm{q}\|_{2,\infty}\leq\alpha} \bigg\{\mathcal{F}(\bm{q}) := \frac{1}{2}\|\bm{\mathrm{div}}\bm{q} + \bb\|^{2} + \frac{1}{2}\|\bb\|^{2} + \frac{\varepsilon}{2\alpha}\|\bm{q}\|^{2}\bigg\}.
    \label{eq:dual_huber_rof}
\end{equation}
It is known that the corresponding primal problem of \eqref{eq:dual_huber_rof} is  the standard Huber-TV denoising \cite{ChambollePock2011} which involves a quadratic smoothing of the $\|\cdot\|_{2,1}$-norm around an $\varepsilon$-neighbourhood of the origin. Among other effects, this reduces the staircasing artifacts of TV, see last image of Figure \ref{fig_shapes_images:main}. 

\begin{figure}[h!]
    \centering
    \begin{subfigure}[b]{8cm}
        \centering
        \includegraphics[width=8cm]{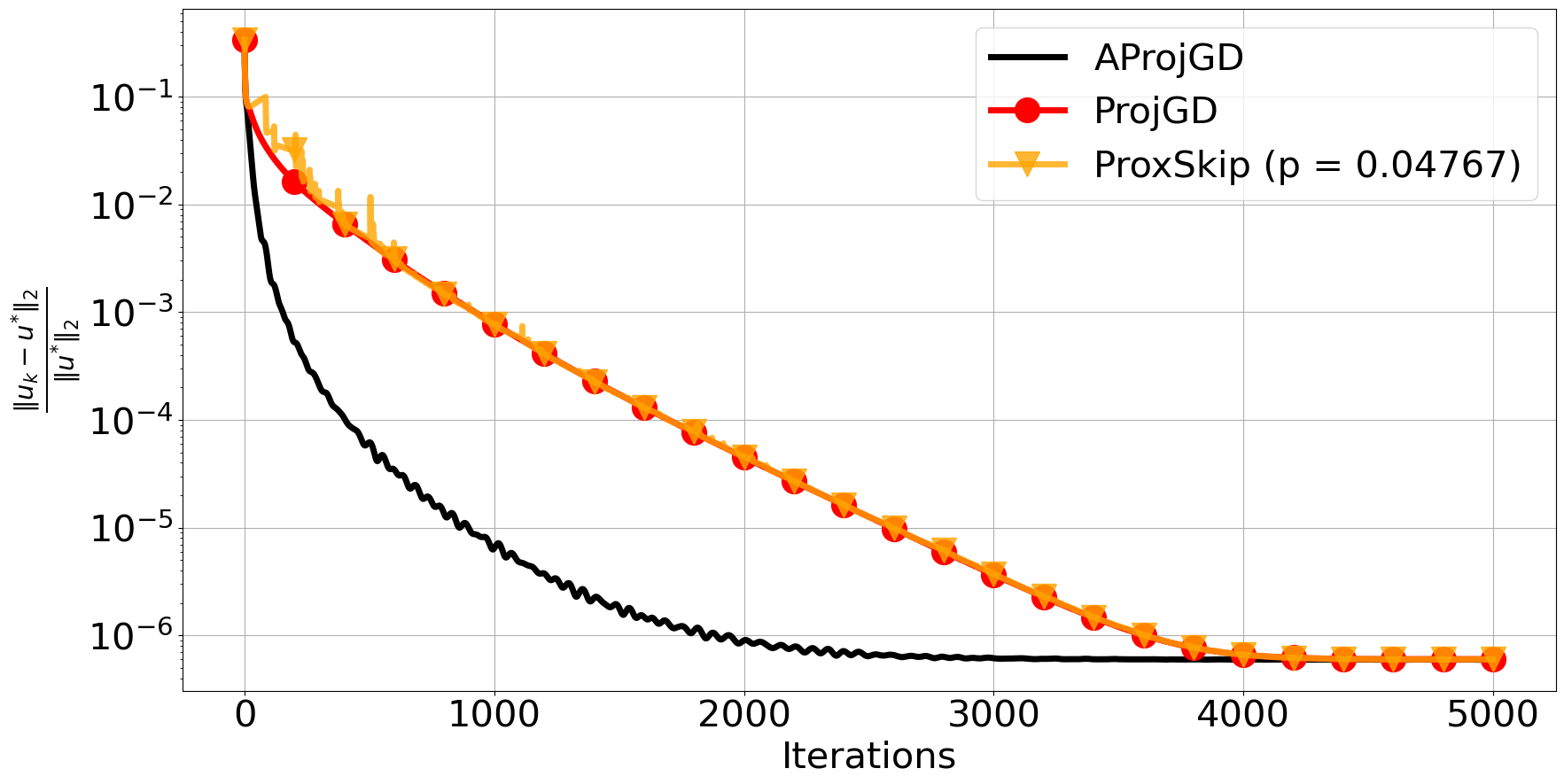}
        \label{fig_shapes_images_strongly_time_iter:sub1}
    \end{subfigure}
    \begin{subfigure}[b]{8cm}
        \centering
        \includegraphics[width=8cm]{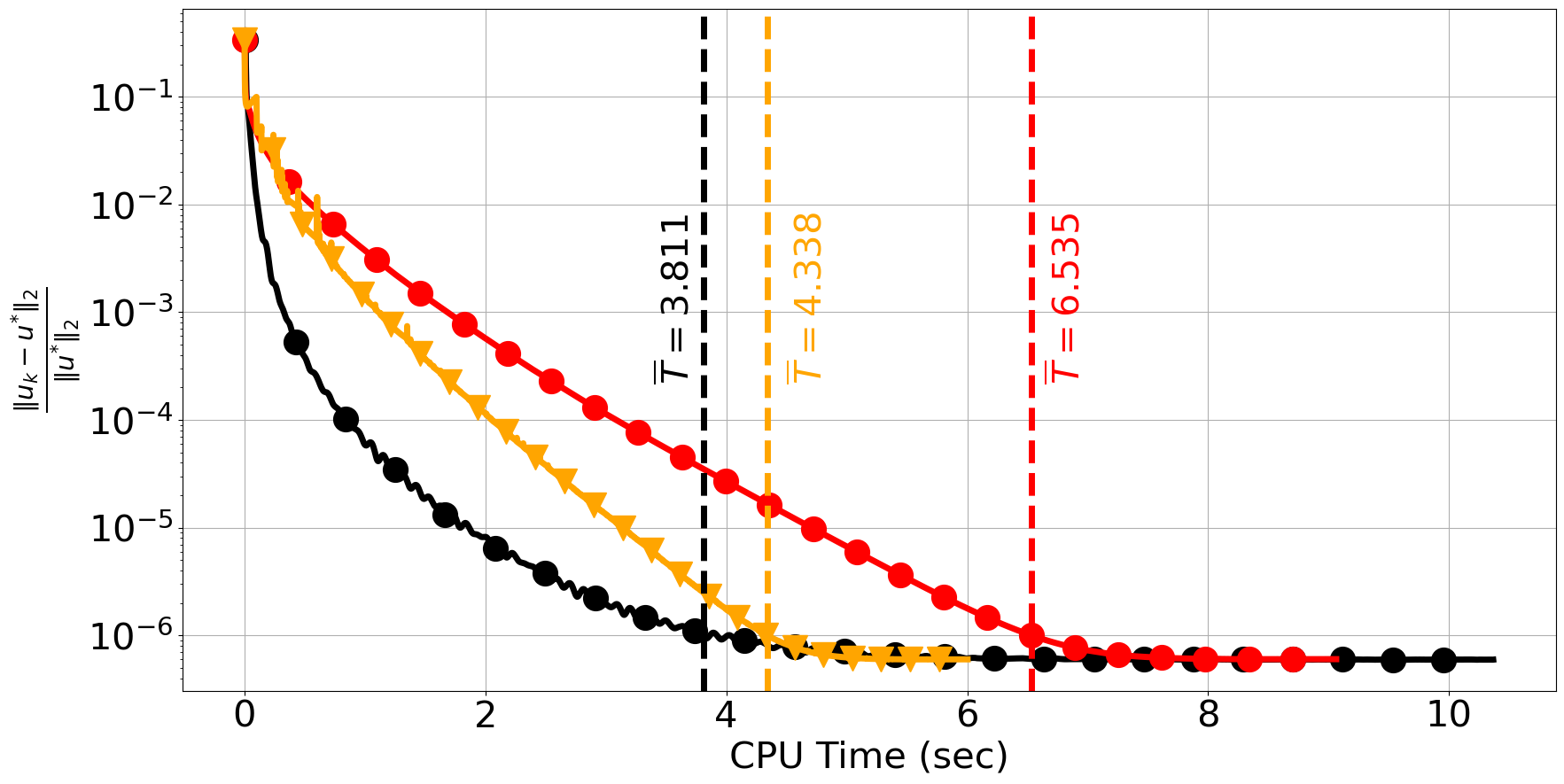}
    \label{fig_shapes_images_strongly_time_iter:sub2}
    \end{subfigure}
    \caption{ProjGD, AProjGD and ProxSkip with optimal $p$ value for the \eqref{eq:dual_huber_rof} problem with respect to iterations (left) and CPU time (right). }    
    \label{fig_shapes_images_strongly_time_iter:main}
\end{figure}

We repeat the same experiment as in the previous section for \eqref{eq:dual_huber_rof} with 5000 iterations and over 30 independent runs of the algorithm, using the probability \ $p=0.04767$, given by \eqref{eq:optimal_prob}.  This results in on average using only 215 projection steps during 5000 iterations. In addition to ProjGD, we also compare its accelerated version, denoted by AProjGD \cite{Nesterov2004} which is essentially Algorithm \ref{alg:GD} with the acceleration step of Algorithm \ref{alg:FISTA}. In Figure \ref{fig_shapes_images_strongly_time_iter:main} (left), we observe that AProjGD performs better  than both ProjGD and ProxSkip with respect to iteration number. While a similar behaviour is observed with respect to the CPU time, see Figure \ref{fig_shapes_images_strongly_time_iter:main} (right), the average time for AProjGD and ProxSkip to reach a relative error less than $10^{-6}$ is nearly identical. One would expect a similar speed-up when AProjGD is combined with skipping techniques but we leave this for future research.

\section{ProxSkip in imaging problems with heavy proximals}\label{sec:heavy}

\subsection{TV deblurring}\label{sec:ProxSkip_TV_deblurring}

We now consider more challenging imaging tasks, involving proximal operators that do not admit closed form solutions, and which thus require computationally intensive iterative solvers. 
\begin{figure}[t!]
    \centering
    \begin{subfigure}[b]{0.24\textwidth}
        \centering
        \includegraphics[width=\textwidth]{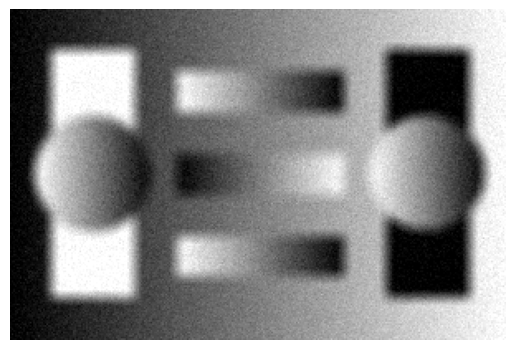}
        \label{fig_shapes_deblurring:sub1}
    \end{subfigure} 
    \begin{subfigure}[b]{0.24\textwidth}
        \centering
        \includegraphics[width=\textwidth]{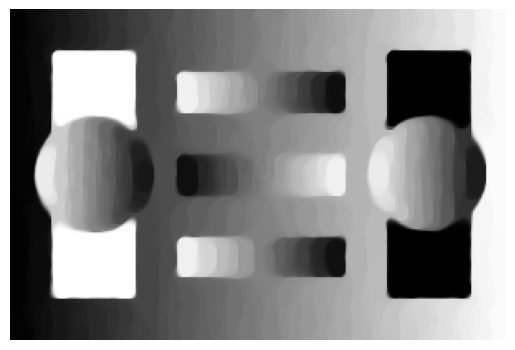}
        \label{fig_shapes_deblurring:sub2}
    \end{subfigure}  
    \begin{subfigure}[b]{0.24\textwidth}
        \centering
        \includegraphics[width=\textwidth]{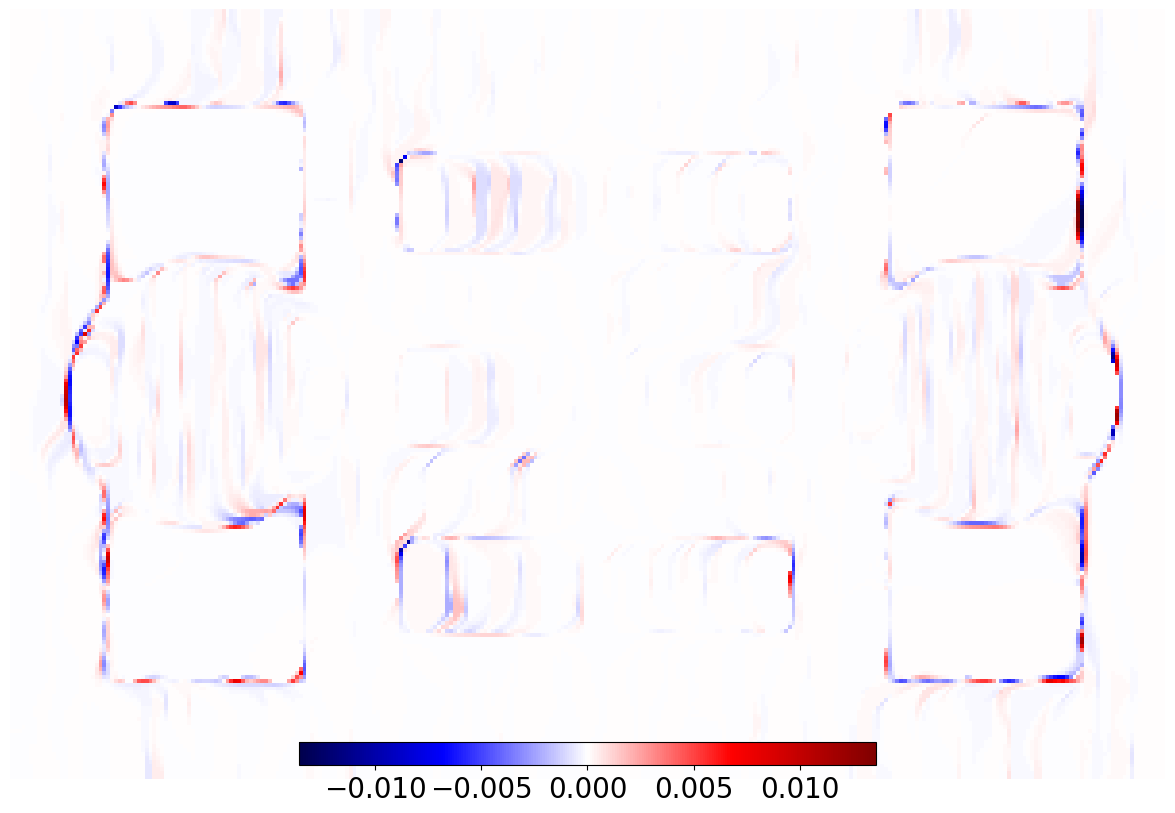}
        \label{fig_shapes_deblurring:sub3}
    \end{subfigure} 
    \begin{subfigure}[b]{0.24\textwidth}
        \centering
        \includegraphics[width=\textwidth]{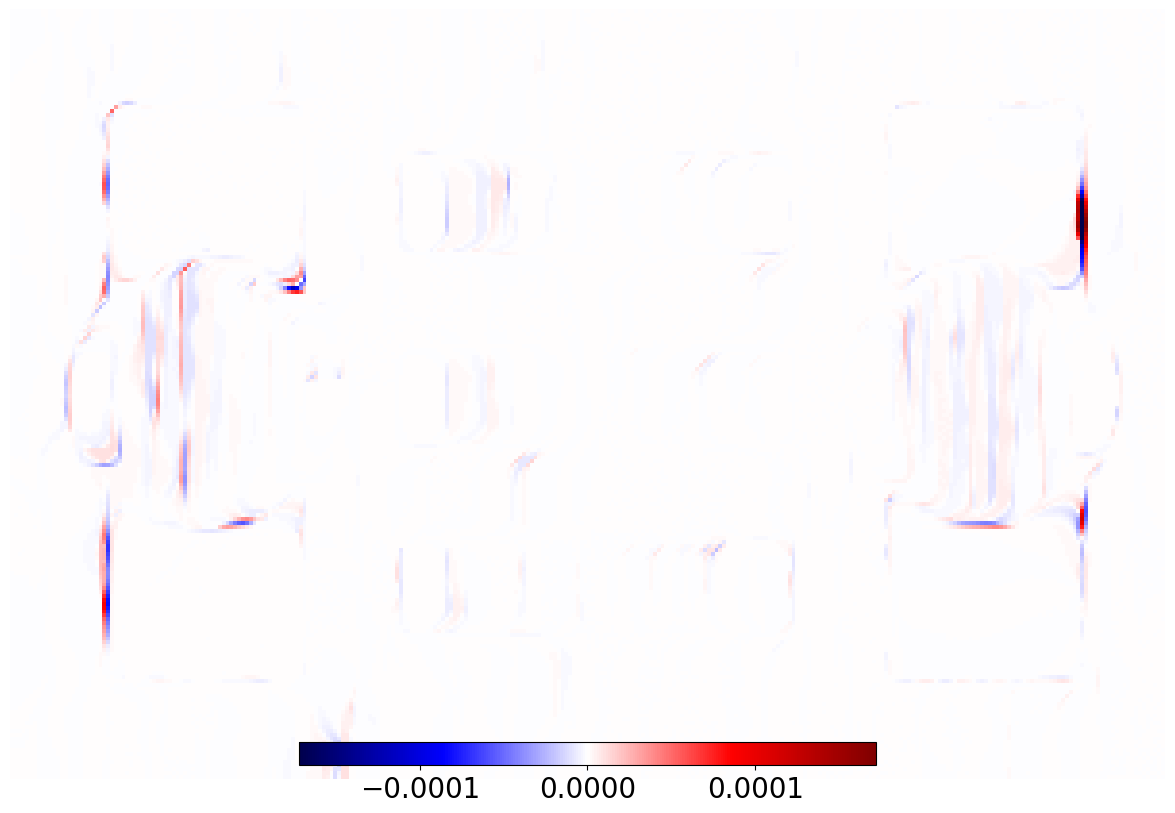}
        \label{fig_shapes_deblurring:sub4}
    \end{subfigure}     
    \caption{Left to right: Noisy and blurry image. TV deblurred  image with $\alpha=0.025$. 
    Difference  $|\bu_{k}-\bu^{\ast}|$ for FISTA when it is less than $\varepsilon=10^{-3}$ and $10^{-5}$. }
    \label{fig_shapes_deblurring:main}
\end{figure}
We start with a deblurring problem in which $\bA$ is a convolution operator, see Figure \ref{fig_shapes_deblurring:main}.

We solve \eqref{mainTV} and compare ISTA, FISTA and  ProxSkip  with different values of $p$. Here the proximal operator corresponds to a TV denoiser for which
we employ AProjGD with a fixed number of iterations
as an inner iterative solver, see next section for other feasible options. In the framework of inexact regularisation another option is to terminate the inner solver based on some metric and predefined threshold  \cite{Rasch2020}. 
As noted therein, the number of required inner iterations typically increases up to $10^3$, as the outer algorithm progresses, leading to higher computational costs over time. To explore both the computationally easy and hard cases, we run the inner solver with 10 and 100 iterations, and use a warm-start strategy. Warm-start is a vital assumption for inexact regularisation as it avoids semi-convergence, where the error stagnates and fails to reach high precision solutions. 
 To avoid biases towards proximal-gradient based solutions,
the ``exact'' solution $\bu^{\ast}$ is computed using 
 200000 iterations of  PDHG  with diagonal preconditioning \cite{ChambollePock_2011_Diagonal}. Outer algorithms are terminated if either 3000 iterations are reached or the relative distance error is less than $\varepsilon=10^{-5}$. The reported CPU time is averaged over 10 runs. 

\begin{figure}[h!]
    \centering
    \begin{subfigure}[t]{8cm}
        \centering
        \includegraphics[width=8cm]{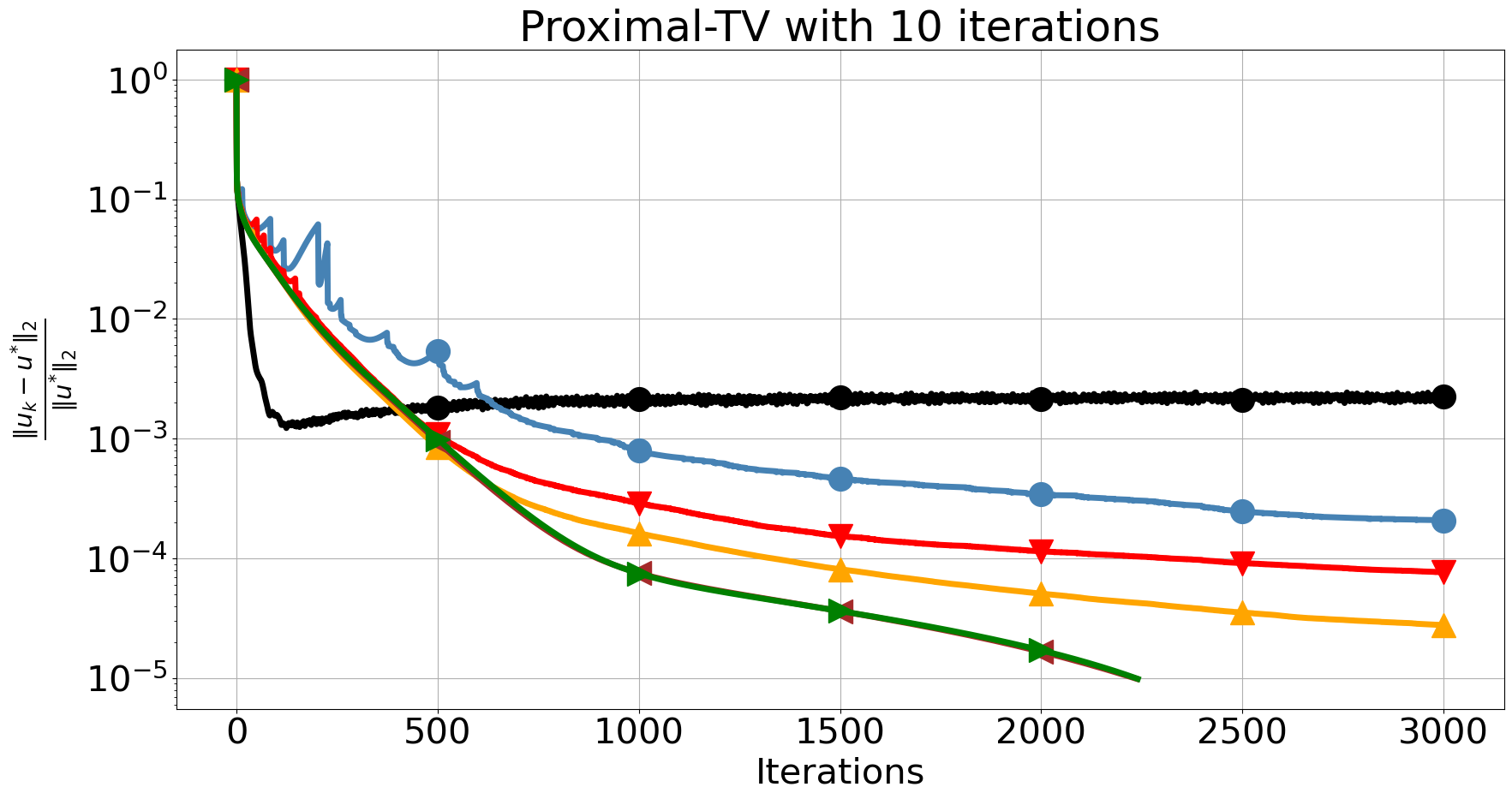}
        \label{fig_shapes_deblurring_comparison:sub1}
    \end{subfigure}
    \begin{subfigure}[t]{8cm}
        \centering
        \includegraphics[width=8cm]{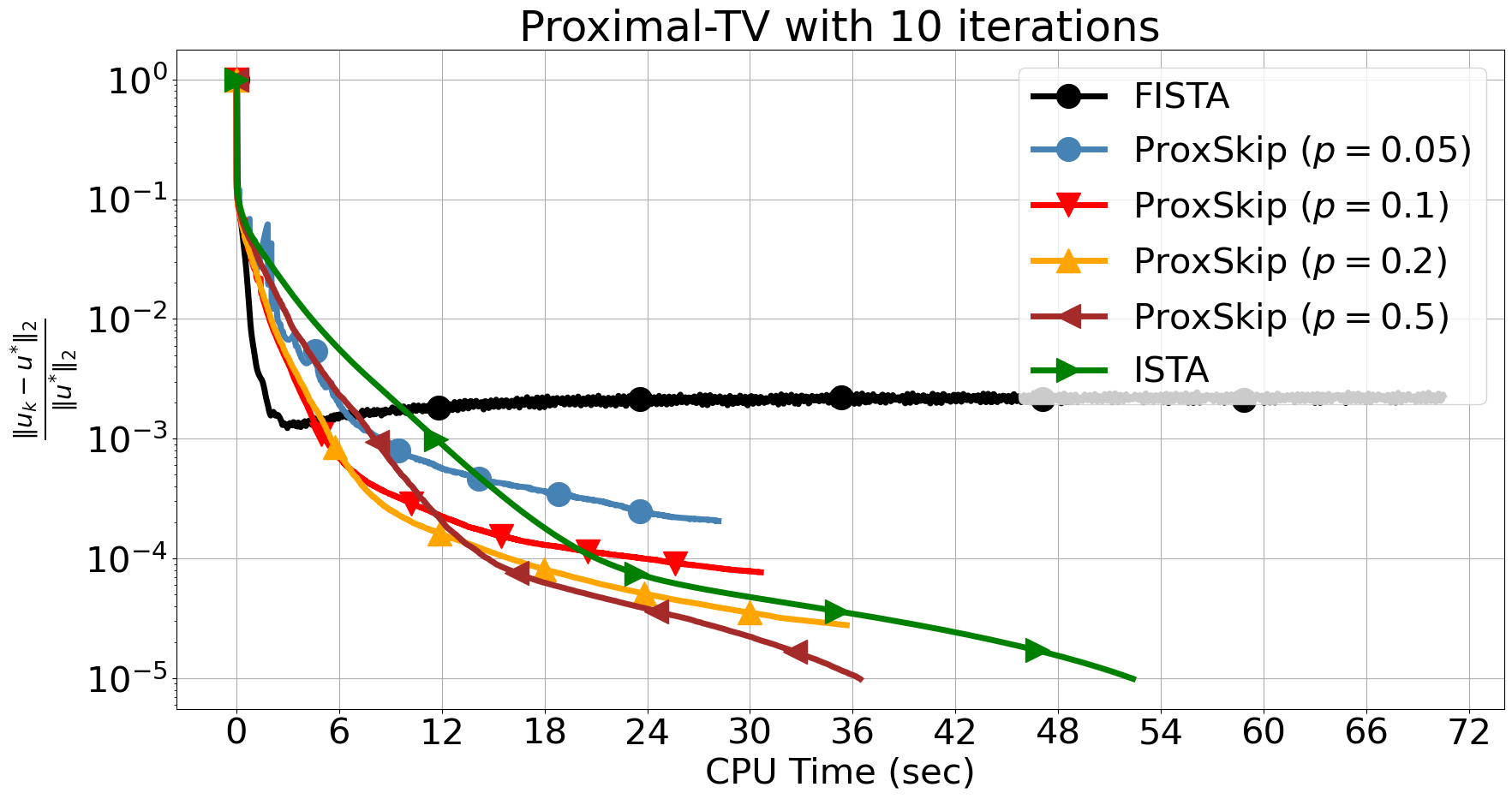}
        \label{fig_shapes_deblurring_comparison:sub2} 
    \end{subfigure}     
    \label{fig:deblurring_iterations}
    \\
    \begin{subfigure}[t]{8cm}
        \centering
        \includegraphics[width=8cm]{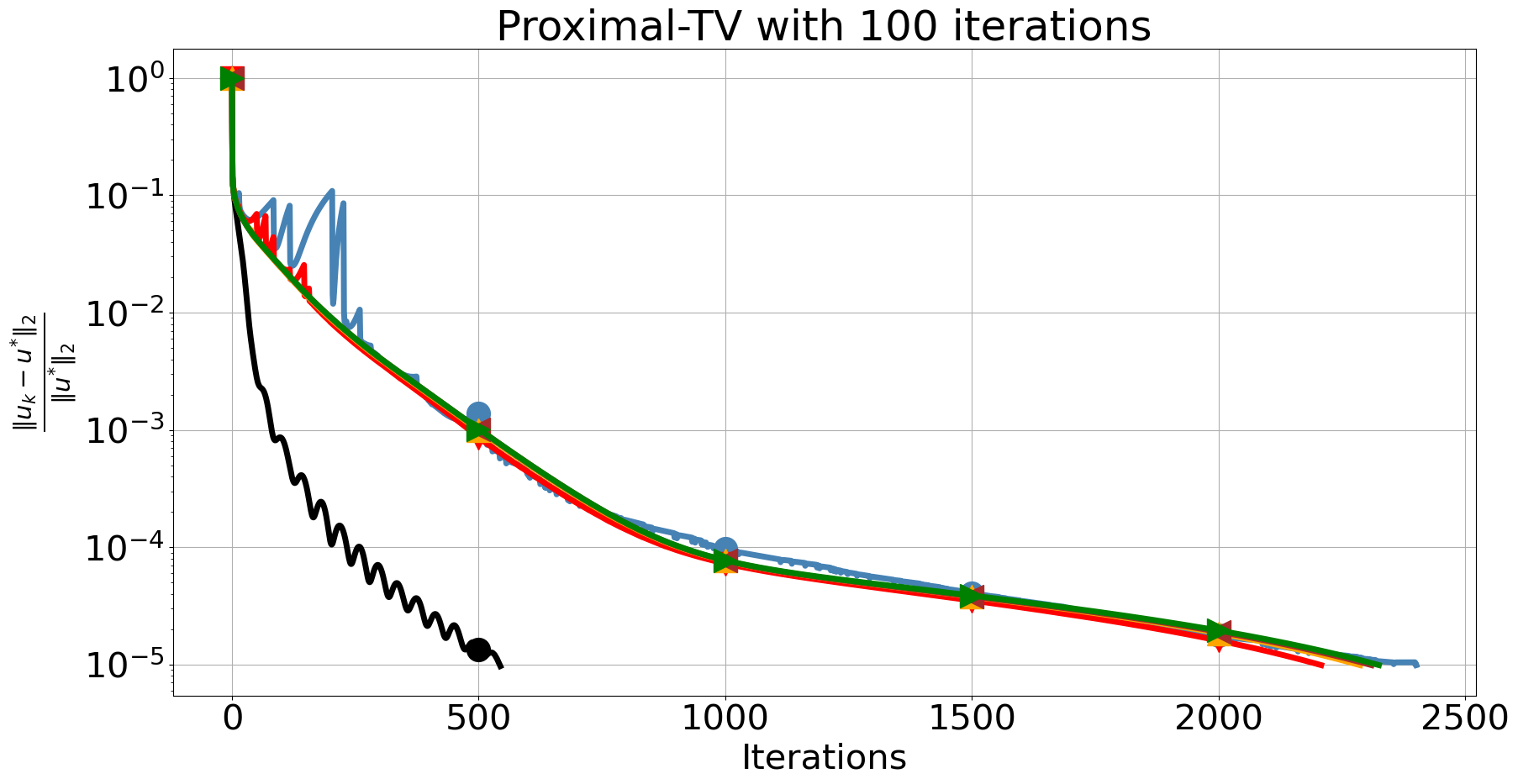}
        \label{fig_shapes_deblurring_comparison:sub5}               
    \end{subfigure}
    \begin{subfigure}[t]{8cm}
        \centering
        \includegraphics[width=8cm]{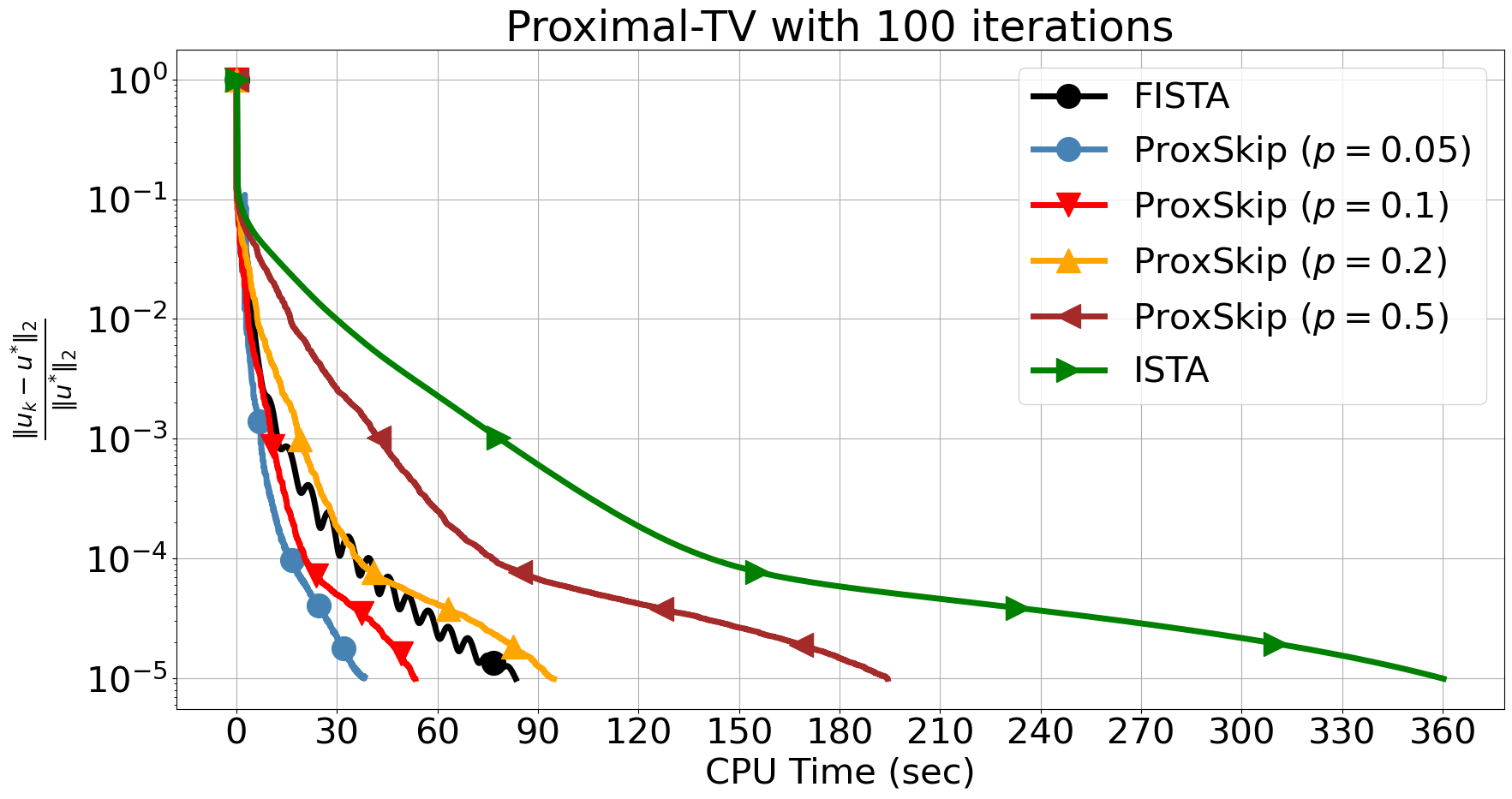}        
        \label{fig_shapes_deblurring_comparison:sub6}  \end{subfigure}  
    \caption{Comparing ISTA, FISTA and Proxskip for multiple values of $p$ for TV deblurring. The proximal of TV is solved using AProjGD using 10 and 100 iterations. ProxSkip outperforms FISTA when $p = 0.05$ and $0.1$.}
    \label{fig_shapes_deblurring_comparison:main}
\end{figure}

In Figure \ref{fig_shapes_deblurring_comparison:main}, we first observe that solving the inner TV problem with 10 iterations seriously affects the convergence of FISTA. Notably, ProxSkip versions are less affected  even though besides the error introduced by the inexact solver, there is also an error from skipping the proximal. By raising the number of inner iterations to  100, and thus increasing the accuracy of the inner solver, we observe that 
FISTA exhibits an early decay albeit with some oscillations and it terminates after around 500 iterations, reaching an accuracy of $\varepsilon=10^{-5}$, see Figure \ref{fig_shapes_deblurring_comparison:main} (bottom-left). On the other hand, ISTA and the ProxSkip require many more iterations to reach the same level of accuracy. 
However, remarkably, in this regime ProxSkip is significantly faster, in terms of CPU time, than  ISTA and even outperforms FISTA when $p=0.05$ and $0.1$, see Figure \ref{fig_shapes_deblurring_comparison:main} (bottom-right). 
In Table \ref{table_high:fig_shapes_deblurring}, we report all the information for the best three algorithms that outperformed FISTA with 100 inner iterations in terms of CPU time.

We note that there are versions of FISTA \cite{ODonoghue2013, Chambolle2015, Liang2022, Aujol2024} with improved performance, also avoiding  oscillations. However, our purpose here is not an exhaustive comparison but to show that a simple version of ProxSkip outperforms the simplest version of FISTA. We anticipate the future development of more sophisticated skip-based algorithms including ones based on accelerated methods.

\begin{table}
\centering
{\renewcommand{\arraystretch}{0.95}
\begin{tabular}{|c|c|c|c|c|c|c|}
\hline
Algorithm & Time $\pm$ Error (sec) & Iterations & \# $\mathrm{prox}_{\mathrm{TV}}$ & Speed-up (\%) \\\hline
ProxSkip -10 ($p=0.5$) & 36.375 $\pm$ 0.893 & 2233 & 1093 & 56.30\\\hline
ProxSkip -100 ($p=0.05$) & 39.43 $\pm$ 0.664 & 2399 & 123 & 21.0\\\hline
ISTA -10 & 52.295 $\pm$ 1.369 & 2237 & 2237 & 37.10 \\\hline
FISTA - 100 & 83.205 $\pm$ 4.933 & 543 & 543 & 0.0 \\\hline
\end{tabular}
}
\caption{The best three algorithms that outperformed FISTA-100 for $\varepsilon=10^{-5}$. 
}
\label{table_high:fig_shapes_deblurring}
\end{table}

\subsection{PDHGSkip}

Skipping the proximal 
  can also be applied to primal-dual type algorithms. Here, for convex $f,g$ and a bounded linear operator $\mathcal{K}$, the  optimisation framework is
\begin{align}
&\min_{\bx\in\mathbb{X}} f(\mathcal{K}\bx) + g(\bx).
\label{eq:pd_objective}
\end{align}
In\cite{Condat2022}, a  skip-version of PDHG \cite{ChambollePock2011}  was proposed, which we denote here by PDHGSkip-1.  It allows not only to skip one of the proximal steps but also the forward or backward operations of $\mathcal{K}$,  depending on the order of the proximal steps, see Algorithm \ref{alg:PDHGSkip_Condat}. In step 4, a Bernoulli operator is used, defined as $\mathcal{B}_{p}(\bx) = \bx/p$ with probability $p$ and $0$ otherwise.
Again,
strong convexity of both $f$ and $g$ is required for convergence. However, our imaging experiments -- in the absence of strong convexity -- revealed a relatively slow performance, even with optimised step sizes $\sigma, \tau$ and $\omega + 1 = 1/p$ according to \cite[Theorem 7]{Condat2022},

\begin{figure}[t!]
\noindent
\begin{minipage}[t]{0.51\textwidth}
\begin{algorithm}[H]
\begin{algorithmic}[1]
\State {\bf Parameters:} $\sigma, \tau, \omega\geq0$, probability\ $p$
\State {\bf Initialize:} $\bx_0 \in \mathbb{X}$, $\by_0 \in \mathbb{Y}$
\For{$k=0,\ldots, K-1$}
\State $\hat{\bx}_{k} = \mathcal{B}_{p}(\mathrm{prox}_{\sigma g}(\bx_{k}-\sigma\mathcal{K}^{T}\by_{k})-\bx_{k})$
\State $\bx_{k+1} = \bx_{k} + \frac{1}{1+\omega}\hat{\bx}_{k}$
\State $\by_{k+1} = \mathrm{prox}_{\tau f^{*}}(\by_{k} + \tau\mathcal{K}(\bx_{k+1} + \hat{\bx}_{k})$
\EndFor  
\end{algorithmic}
\caption{PDHGSkip-1 \cite{Condat2022}}
\label{alg:PDHGSkip_Condat}
\end{algorithm}    
\end{minipage}\hfill
\begin{minipage}[t]{0.49\textwidth}
\begin{algorithm}[H]
\begin{algorithmic}[1]
\State {\bf Parameters:} $\sigma, \tau >0$, prob.\ $p$ 
\State {\bf Initialize:} $\bx_0, \bh_{0}  \in \mathbb{X}$, $\by_0\in \mathbb{Y}$
\For{$k=0,\ldots, K-1$}
\State $\hat{\bx}_{k} = \bx_{k} - \sigma(\mathcal{K}^{T}\by_{k} - \bh_{k})$
\If{$\theta_k=1$} 
\State  $\bx_{k+1} = \mathrm{prox}_{\frac{\sigma}{p}g}\Big(\hat \bx_{k} - \frac{\sigma}{p}{\bh_k} \Big)$ 
\Else
$\;\;\bx_{k+1} = \hat \bx_{k}$ \hfill 
\EndIf
\State $\Bar{\bx}_{k+1} =2\bx_{k+1} -\bx_k$
\State $\by_{k+1} = \mathrm{prox}_{\tau f^*}\left(\by_k+\tau \mathcal{K}\Bar{\bx}_{k}\right)$
\State ${ \bh_{k+1}} = {\bh_k} + \frac{p}{\sigma}(\bx_{k+1} - \hat \bx_{k+1})$ 
\EndFor  
\end{algorithmic}
\caption{PDHGSkip-2 (ours)}
\label{alg:PDHGSkip}
\end{algorithm}    
\end{minipage}
\end{figure}

\begin{figure}[h!]
    \centering
    \begin{subfigure}[t]{8cm}
        \centering
        \includegraphics[width=8cm]{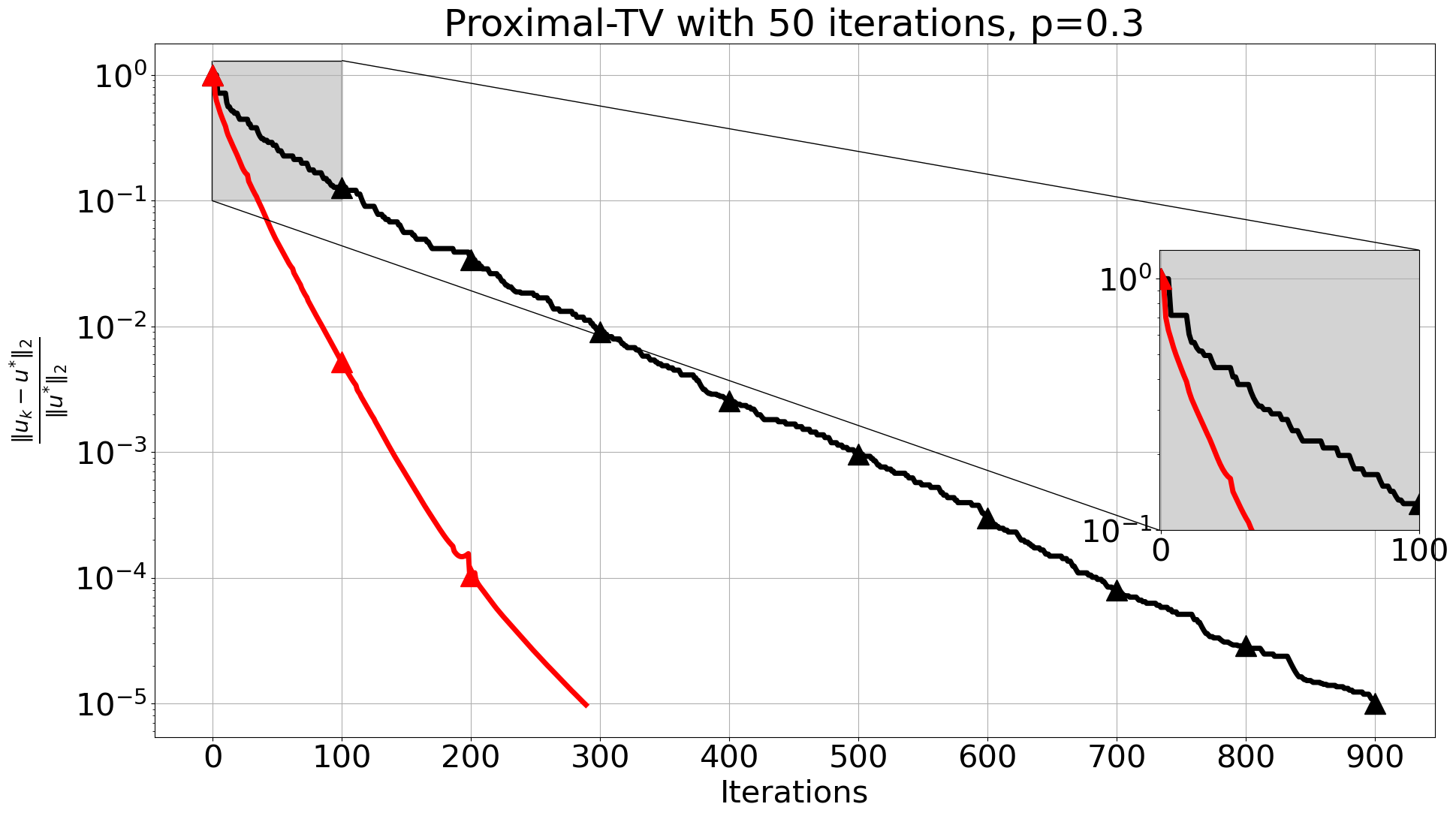}
    \end{subfigure}
    \begin{subfigure}[t]{8cm}
        \centering
        \includegraphics[width=8cm]{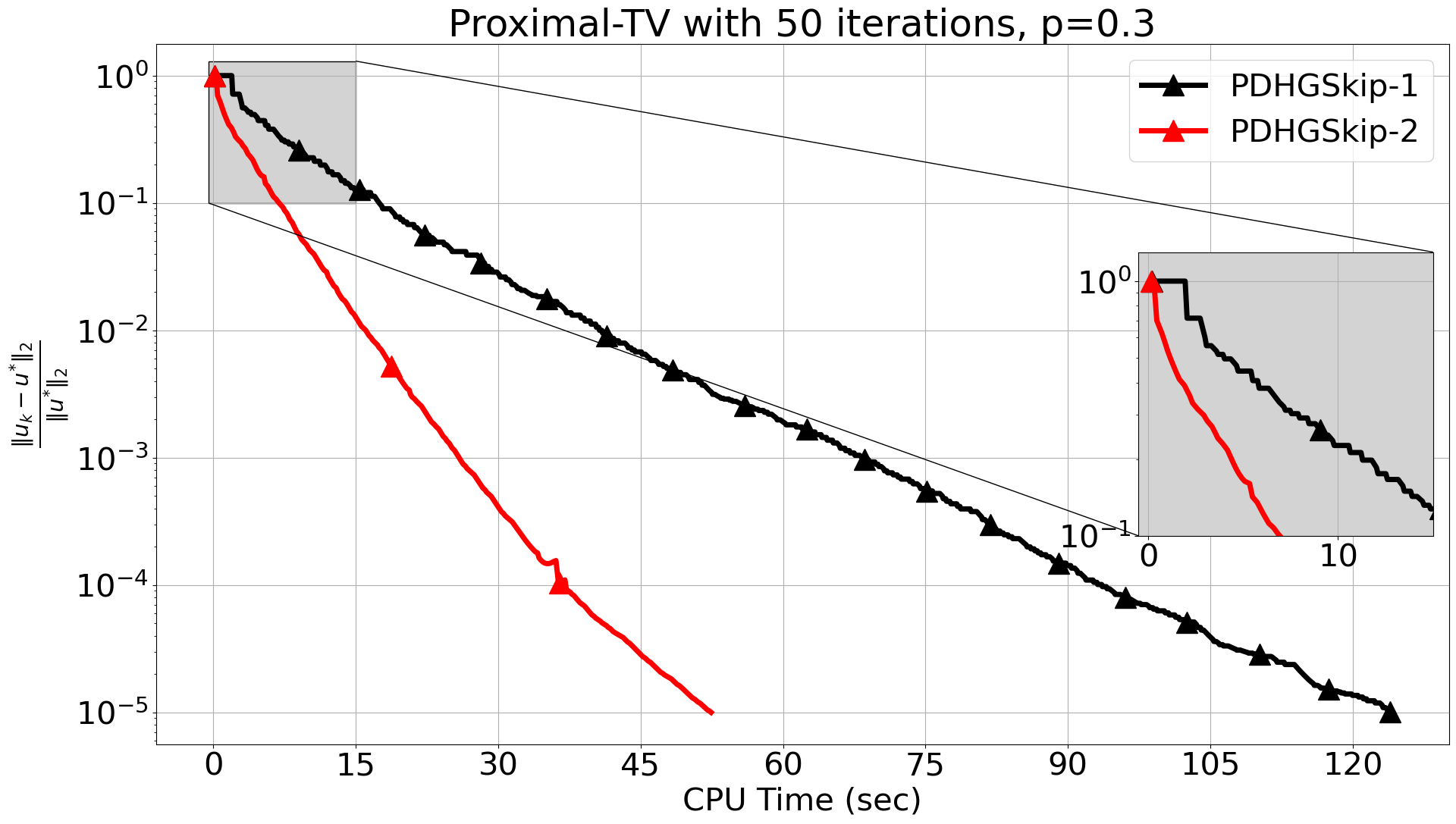}
    \end{subfigure}            
    \caption{Comparing PDHGSkip-1 by \cite{Condat2022} and our proposed PDHGSkip-2 for $p=0.3$ for the tomography problem of Section \ref{sec:ProxSkip_TV_tomography}, using  50 iterations for the inner  solver for the proximal (a TV denoising problem).}
\label{fig_finden_tomography_condat_pdhgskip_comparison:main}  
\end{figure}

To address the slow convergence, we introduce a modification: PDHGSkip-2, see Algorithm \ref{alg:PDHGSkip}. The difference to PDHGSkip-1 is that the adjoint $\mathcal K^T$ and the proximal operator of $g$ are now separated, see steps 4-7. Note that in both versions, when $p=1$ the control variable $\bh$ vanishes and PDHG is recovered. To illustrate the difference of the two versions in practice, we show in Figure \ref{fig_finden_tomography_condat_pdhgskip_comparison:main} a simple comparison on tomography, see next section for details. Apart from the clear acceleration of PDHGSkip-2 over PDHGSkip-1. We also observe a staircasing pattern for the relative error for PDGHSkip-1, see detailed zoom of the first iterations. This  is expected since in most iterations, where the proximal step is skipped, one variable vanishes without contributing to the next iterate. Hence, the update remains unchanged.

In general,  $\ell^{2}$-$\mathrm{TV}$  problems can be solved using  implicit or explicit formulations of PDHG. In the implicit case, $f$ is the $\|\cdot\|_{2}^{2}$ term, $\mathcal{K}=\bA$ and $g$ is the $\mathrm{TV}$ term. In the explicit case, $f$ is a separable sum of $\|\cdot\|_{2}^{2}$ and $\|\cdot\|_{2,1}$ composed with block operator $\mathcal{K} = [\bA, \dmat]^{T}$ and $g$ can be a zero function or a non-negativity constraint. Here, every proximal step (Steps 6, 9 in Algorithm \ref{alg:PDHGSkip}) has an analytic solution. This significantly reduces the cost per-iteration, but also requires  more iterations to reach a desired accuracy \cite{Burger2016}. Inexact regularisation is usually preferred and typically reduces the number of iterations. Hence, the number of calls of the forward and backward operations of $\bA$ is reduced, which in certain applications gives a considerable speed-up.

\subsection{TV Tomography reconstruction with real-world data}
\label{sec:ProxSkip_TV_tomography}

\begin{figure}[t!]
    \centering
    \begin{subfigure}[b]{0.25\textwidth}
        \centering
        \includegraphics[width=\textwidth]{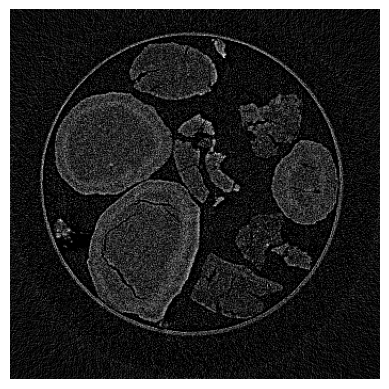}
        \label{fig_tomography_images:sub1}
    \end{subfigure} 
    \begin{subfigure}[b]{0.25\textwidth}
        \centering
        \includegraphics[width=\textwidth]{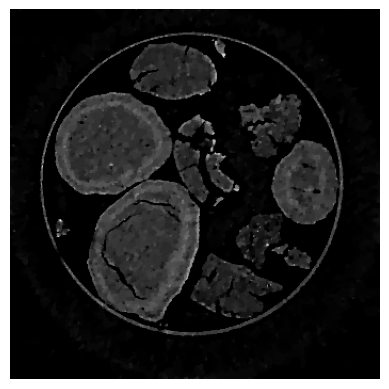}
        \label{fig_tomography_images:sub2}
    \end{subfigure}     
    \caption{FBP (left) and TV (right) reconstruction using diagonal preconditioned PDHG for 200000 iterations. Regularisation parameter is manually set to balance noise reduction with feature preservation.}
    \label{fig_tomography_images:main}
\end{figure}

For our final case study, we solve  \eqref{mainTV} under a non-negativity constraint for $\bu$ for a real-world tomographic reconstruction. Here $\bA$ is the discrete  Radon transform and $\bb$ is a noisy sinogram of a real chemical imaging tomography dataset, representing post partial oxidation of methane  reaction Ni-Pd/CeO2-ZrO2/Al2O3 catalyst \cite{Matras2021, Vamvakeros2018}, see Figure \ref{fig_tomography_images:main}. The initial dataset was acquired for 800 projection angles with 695$\times$695 detector size and 700 vertical slices. For demonstration purposes and to be able to perform multiple runs for computing more representative CPU times, we consider one vertical slice with half the projections and $2\times$ rebinned detector size. Same conditions for the inexact solver and stopping rules are used as in previous section. For these optimisation problems, one can use algorithms that fit both the general  frameworks \eqref{eq:fbs_objective} and \eqref{eq:pd_objective}, see \cite{Anthoine2012}.

\begin{figure}[h!]
    \centering
    \begin{subfigure}[t]{8cm}
        \centering
        \includegraphics[width=8cm]{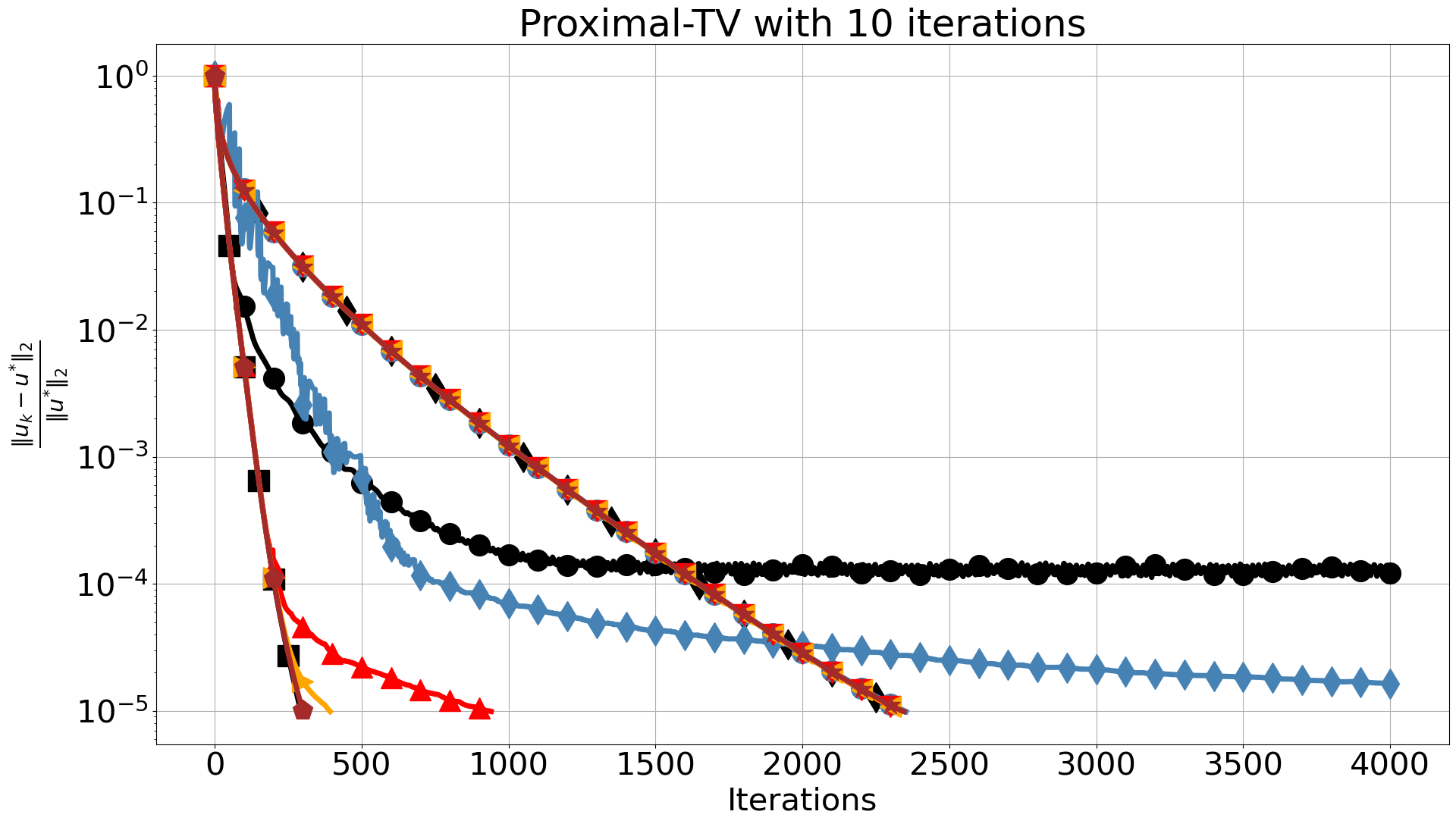}
        \label{fig_finden_tomography_comparison:sub1}
    \end{subfigure}
    \begin{subfigure}[t]{8cm}
        \centering
        \includegraphics[width=8cm]{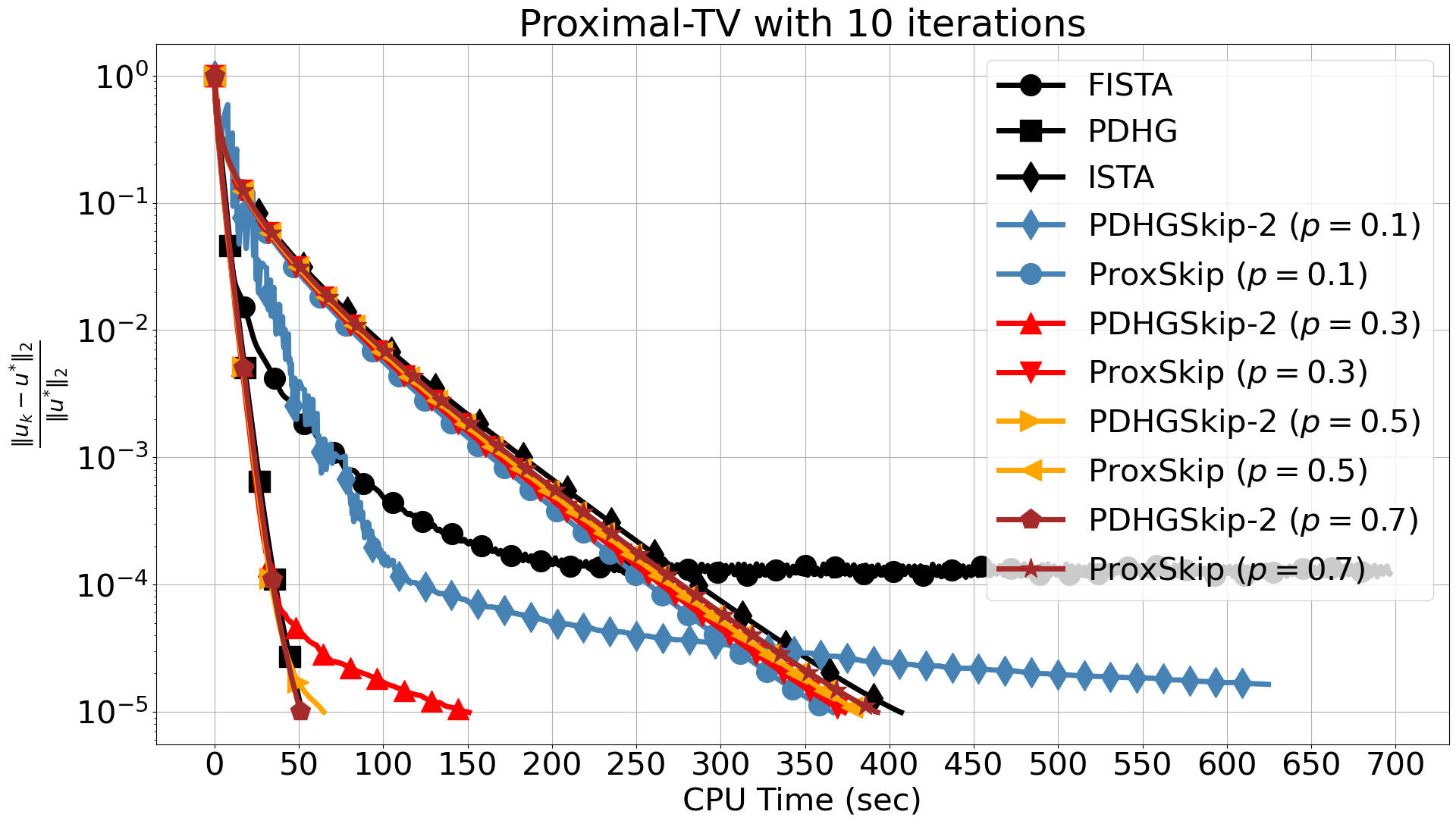}
        \label{fig_finden_tomography_comparison:sub2} 
    \end{subfigure}     
    \label{fig:deblurring_iterations}
    \\
    \begin{subfigure}[t]{8cm}
        \centering
        \includegraphics[width=8cm]{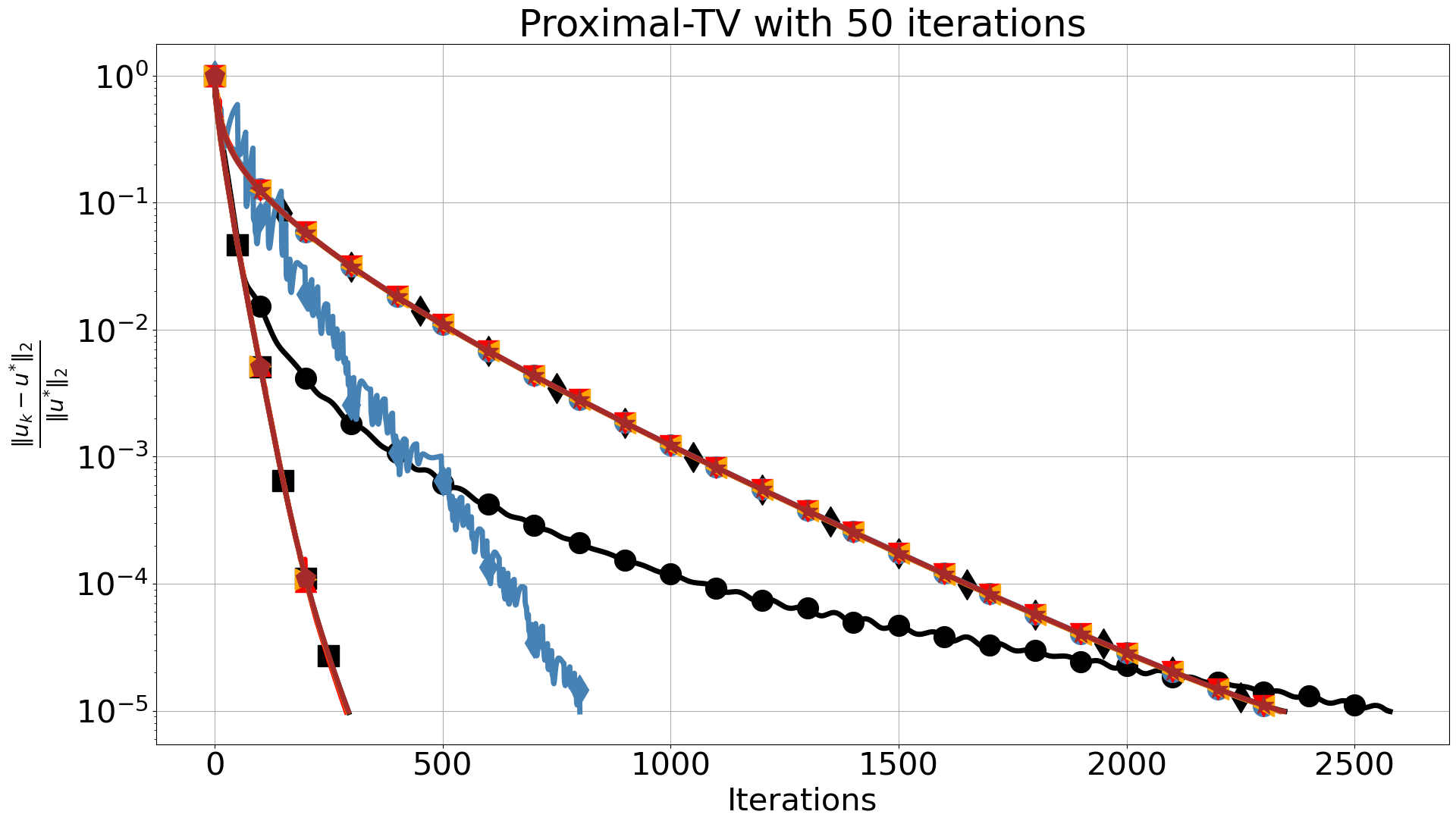}
        \label{fig_finden_tomography_comparison:sub3}
    \end{subfigure}
    \begin{subfigure}[t]{8cm}
        \centering
        \includegraphics[width=8cm]{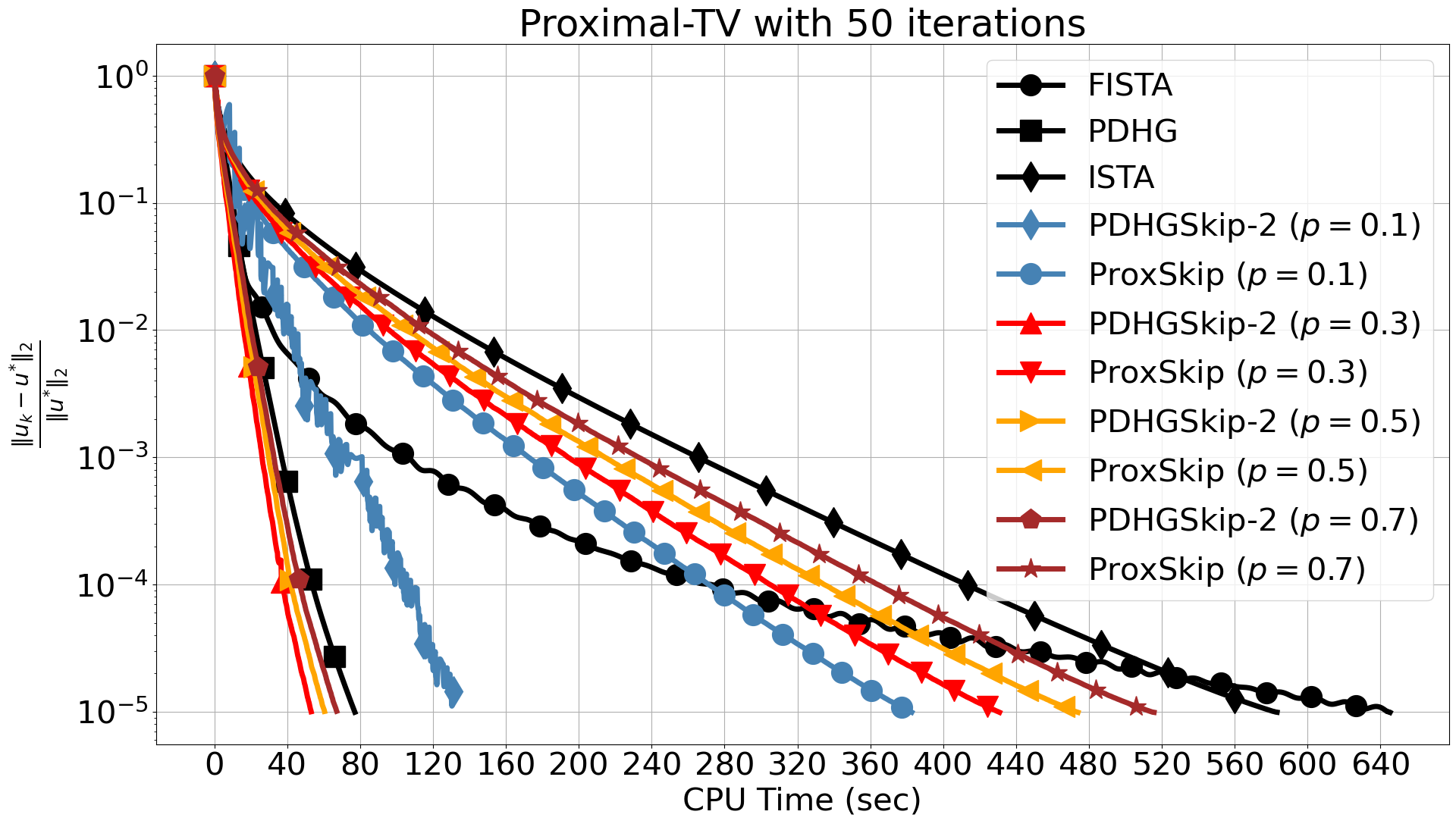}
        \label{fig_finden_tomography_comparison:sub4}
    \end{subfigure}  
    \\
    \begin{subfigure}[t]{8cm}
        \centering
        \includegraphics[width=8cm]{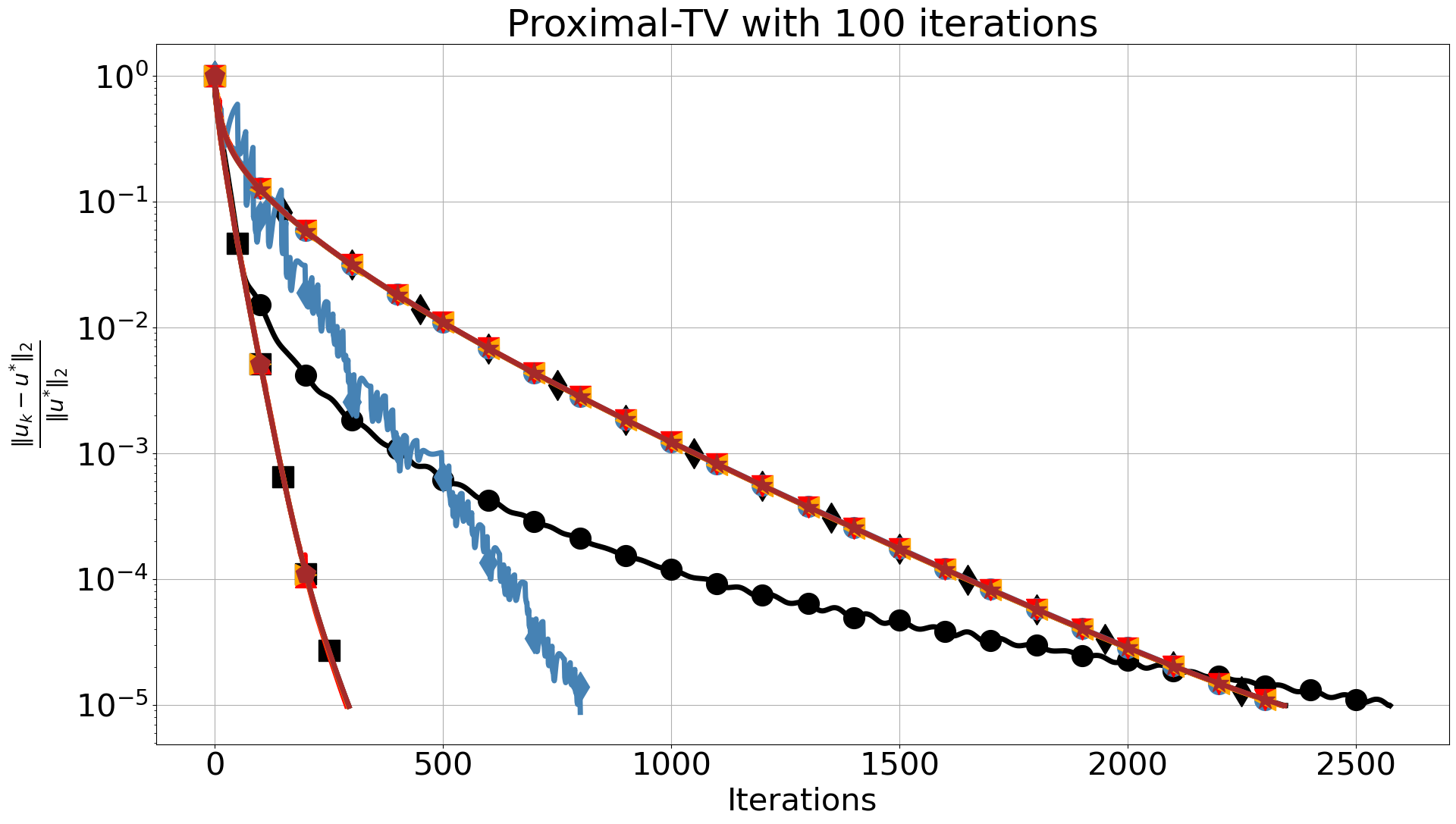}
        \label{fig_finden_tomography_comparison:sub5}
    \end{subfigure}
    \begin{subfigure}[t]{8cm}
        \centering
        \includegraphics[width=8cm]{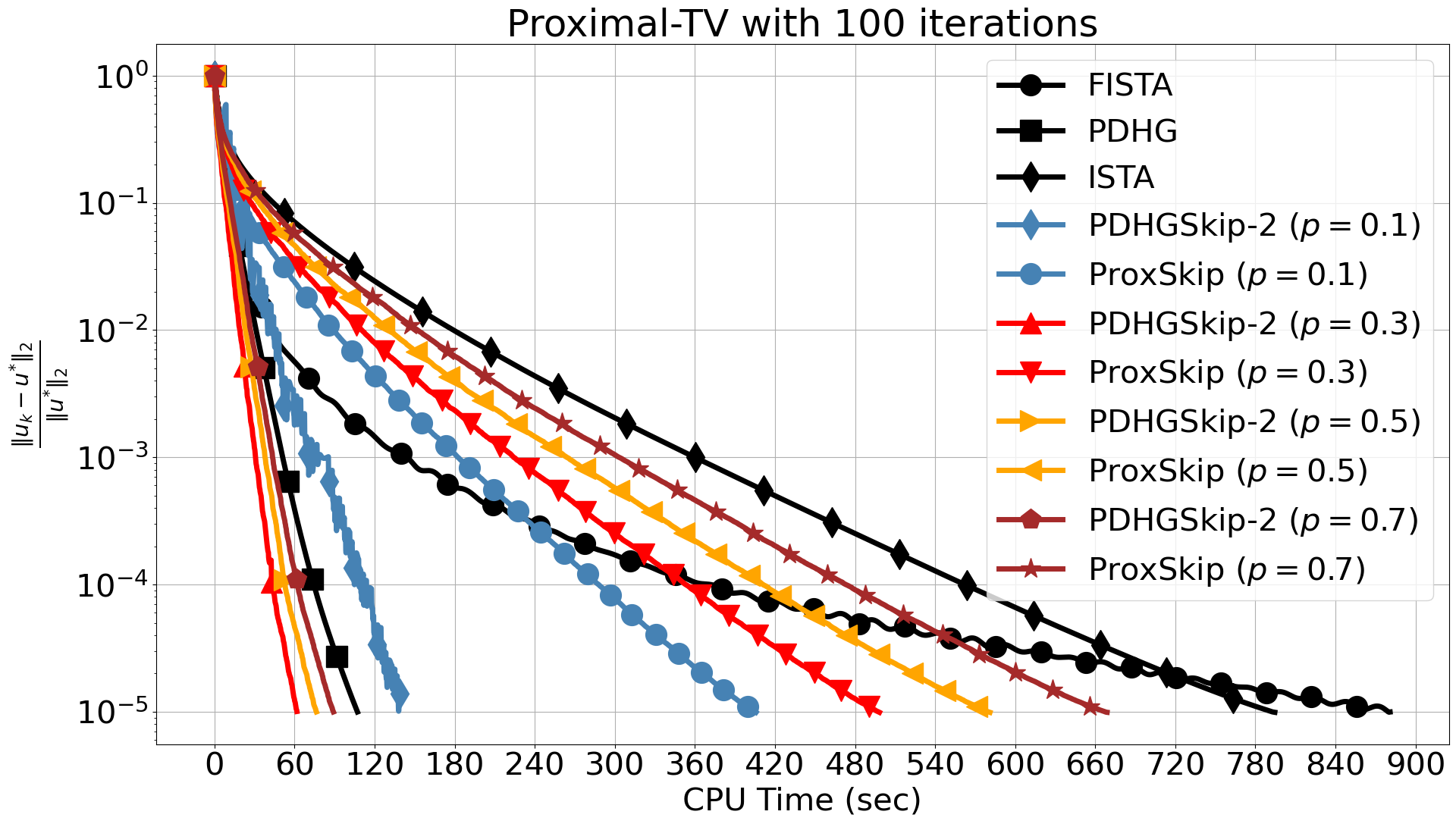}
        \label{fig_finden_tomography_comparison:sub6}
    \end{subfigure}  
    \caption{Comparing ISTA, FISTA, ProxSkip and PDHGSkip-2 for multiple values of $p$ for the TV tomography  problem. The proximal of TV is solved using AProjGD and for 10, 50 and 100 iterations.}
    \label{fig_finden_tomography_comparison:main}
\end{figure}

The two algorithms use the same values for step sizes $\sigma$ and $\tau$, satisfying $\sigma\tau\|\mathcal{K}\|^{2}\leq1$. 
In the comparisons presented in Figure \ref{fig_finden_tomography_comparison:main}, we observe a similar trend for the proximal-gradient based algorithms as in Section \ref{sec:ProxSkip_TV_deblurring}. 
When we use 10 iterations for the inner solver,
FISTA fails to reach the required accuracy and stagnates, which is not the case for the ProxSkip algorithms even for the smallest $p$. 
This is accompanied by a significant CPU time speed-up, which further increases when more inner iterations are used. 
In fact, by increasing the number of inner iterations, the computational gain is evident with around 90\% speed-up compared to FISTA. 

The best overall performance with respect to CPU time is achieved by PDHGSkip-2.
There, distance errors are identical to the case $p=1$, in terms of iterations, except for the extreme case $p=0.1$ which oscillates in the early iterations. For 10 inner iterations and $p=0.3, 0.5$ we observe a slight delay towards $\varepsilon=10^{-5}$, due to error accumulation caused by skipping the proximal operator and limited accuracy of the inner solver. This  is corrected when we increase the number of inner iterations. Moreover, we see no difference with respect to CPU time, for 10 inner iterations and $p=0.7, 1.0$. This is expected since the computational cost to run 10 iterations of AProjGD is relatively low. However, it demonstrates that we can have the same reconstruction using the proximal operator 70\% of the time. Finally, the computational gain is more apparent when we increase the number of inner iterations hence the computational cost of the inner solver. 

We note that such expensive inner steps are  used by   open source imaging libraries and are solved with different algorithms and stopping rules.
In CIL \cite{Jorgensen2021,Papoutsellis2021}, the AProjGD is used to solve \eqref{eq:dual_TV} with 100 iterations as default stopping criterion. In PyHST2 \cite{Mirone2014}, AProjGD is used with 200 iterations and a duality gap 
is evaluated every 3 iterations. In TIGRE  \cite{Biguri2016}, the PDHG algorithm with adaptive step sizes is applied to \eqref{eq:dual_TV}, \cite{Zhu2008},
with 50 iterations. In DeepInv
\footnote{https://github.com/deepinv/deepinv}
 PDHG is applied to \eqref{mainTV} (with $\bA = \mathbf{Id}$) with 1000 iterations and error distance between two consecutive iterates. 
Finally, in Tomopy\cite{Grsoy2014}, PDHG is applied to \eqref{mainTV} (with $\bA  = \mathbf{Id}$) and the number of iterations is specified by the user.
All these default options are optimised and tested for particular real-world tomography applications like the ones encountered in synchroton facilities for instance.
Alternatively, one can avoid inexact solvers for TV denoising \cite{Kamilov2016}. There, the proximal is replaced by a combination of wavelet and scaling transforms and is computed using a componentwise shrinkage operator. Even in this case, we expect an improvement by skipping this operator. Overall, we anticipate a computational gain proportional to the cost savings achieved by omitting specific mathematical operations.

\section{Code Reproducibility}
The code and datasets need to reproduce the results will be made available upon acceptance of this paper. Also, we would like to highlight that all the experiments were tested in three different computing platforms under different operating systems. For this paper, we use an Apple M2 Pro, 16Gb without GPU use to avoid measuring data transferring time between the host and the device which can be misleading.

\section{Discussion and Future work}

In this paper, we explored the use of the ProxSkip algorithm for imaging inverse problems. This algorithm allows to skip costly proximal operators that are usually related to the regulariser without impacting the convergence and the final solution. In addition, we presented a new skipped version of PDHG, a more flexible algorithm, which can be useful when L-smoothness assumption is not satisfied, e.g., Kullback-Leibler divergence and its convergence is left for future work. Although, we demonstrated that avoiding computing the proximal leads to better computational times, this speed-up can be further increased when dealing with larger datasets and more costly regularisers, e.g. TGV. 
Additionally to the skipping concept, one can combine stochastic optimisation methods and different variance reduced estimators that use only a subset of the data per iteration. For example, in tomography applications, where the cost per iteration is mostly dominated by the forward and backward operations, one can  randomly select a smaller subset of projection angles in addition to a random evaluation of the proximal operator. In this scenario, the cost per iteration is significantly decreased and from ongoing experiments outperformed deterministic algorithms in terms of CPU time. Finally, we note that we could achieve a further computational gain  if we use the skipping concept to the inner solver as well. Notice for example that  Algorithm \eqref{alg:PDHGSkip} is designed by default to skip the proximal operator of $g$. For the TV denoising problem, the projection step is avoided if 
the dual formulation \eqref{eq:dual_TV} is used, as presented in Section \ref{sec:dualTV}. On the other hand, the proximal operator, related to the fidelity term is avoided if we use the primal formulation.\\
Overall, proximal based algorithms presented here and possible future extensions open the door to revisiting a range of optimisation solutions for a plethora of imaging modalities, that were developed over the last decades, with the potential to greatly reduce actual computation times.

\subsubsection*{Acknowledgements}
E. P. acknowledges funding through the Innovate UK Analysis for Innovators (A4i) program: \emph{Denoising of chemical imaging and tomography data} under which the experiments were initially conducted. E. P. acknowledges also CCPi (EPSRC grant EP/T026677/1). \v{Z}. K. is supported by the UK EPSRC grant EP/X010740/1.

\bibliographystyle{ieeetr}
\bibliography{arxiv_version}

\end{document}